\documentclass[conf]{new-aiaa}
\usepackage[utf8]{inputenc}
\usepackage{bm}
\usepackage{color}
\usepackage{graphicx}
\usepackage{amsmath}
\usepackage[version=4]{mhchem}
\usepackage{siunitx}
\usepackage{subfigure}
\usepackage{longtable,tabularx}
\usepackage[table]{xcolor}
\usepackage{booktabs}
\usepackage[normalem]{ulem}

\newcounter{lnote}

\usepackage[ruled]{algorithm2e}

\SetKwProg{Fn}{}{}{end} 

\newcommand{\ii}{\mathbf{i}}
\newcommand{\jj}{\mathbf{j}}
\newcommand{\ee}{\mathbf{e}}
\title{Uncertainty Quantification of Ship Resistance via Multi-Index Stochastic Collocation and Radial Basis Function Surrogates: \\
A Comparison}

\author{
Chiara Piazzola\footnote{Postdoctoral Research Fellow, CNR-IMATI Pavia, Via Adolfo Ferrata 5/A, 27100} and 
Lorenzo Tamellini\footnote{Research Scientist, CNR-IMATI Pavia, Via Adolfo Ferrata 5/A, 27100}}
\affil{CNR-IMATI, National Research Council-Institute for Applied Mathematics and Information Technologies  \\ ``E. Magenes'', Pavia, Italy}
\author{
Riccardo Pellegrini\footnote{Postdoctoral Research Fellow, CNR-INM Rome, Via di Vallerano 139, 00128}, 
Riccardo Broglia\footnote{Research Scientist, CNR-INM Rome, Via di Vallerano 139, 00128}, 
Andrea Serani\footnote{Research Scientist, CNR-INM Rome, Via di Vallerano 139, 00128}, and
Matteo Diez\footnote{Senior Research Scientist, CNR-INM Rome, Via di Vallerano 139, 00128, AIAA Member, Email: matteo.diez@cnr.it}}
\affil{CNR-INM, National Research Council-Institute of Marine Engineering, Rome, Italy}

\begin{document}

\maketitle

\begin{abstract}
This paper presents a comparison of two methods for the forward uncertainty quantification (UQ) of complex industrial problems. 
Specifically, the performance of Multi-Index Stochastic Collocation (MISC) and adaptive multi-fidelity Stochastic Radial Basis Functions (SRBF) surrogates is assessed for the UQ of a roll-on/roll-off passengers ferry advancing in calm water and subject to two operational uncertainties, namely the ship speed and draught.
The estimation of expected value, standard deviation, and probability density function of the (model-scale) resistance is presented and discussed; the required simulations are obtained by the in-house unsteady multi-grid Reynolds Averaged Navier-Stokes (RANS)
{solver $\chi$navis}. 
Both MISC and SRBF use as multi-fidelity levels the evaluations on the different grid levels intrinsically employed by the RANS solver for multi-grid acceleration; four grid levels are used here, obtained as isotropic coarsening of the initial finest mesh.
The results suggest that MISC could be preferred when only limited data sets are available. For larger data sets both MISC and SRBF represent a valid option, with a slight preference for SRBF, due to its robustness to noise.
\end{abstract}

%%%%%%%%%%%%%%%%%%%%%%%%%%%%%%%%%%%%%%%%%%%%%%%%%%%%%%%%%%%%%%%%%%%%%%%%%%

\section{Introduction}
\lettrine{S}{hip} performance depends on design and operational/environmental parameters. The accurate prediction of significant design metrics (such as resistance and powering requirements; seakeeping, maneuverability, and dynamic stability; structural response and failure) requires prime-principles-based high-fidelity computational tools (e.g., computational fluid/structural dynamics, CFD/CSD), especially for innovative configurations and off-design conditions. These tools are generally computationally expensive, making the exploration of {the spaces of design and operational parameters (as done e.g. in optimization and uncertainty quantification, UQ, respectively)} a technological challenge.

The development and application of UQ methods for vehicle problems (including ships) were the subject of the NATO Science and Technology Organization, Applied Vehicle Technology group  AVT-191 ``Application of Sensitivity Analysis and Uncertainty Quantification to Military Vehicle Design,'' where the UQ of a high-speed catamaran in irregular head waves was performed via both CFD computations \cite{he2013urans,diez2018statistical} and towing-tank experiments \cite{durante2020accurate}. UQ methods for ship operational parameters are an essential element of reliability-based and robust design optimization for vessels sailing in real-world stochastic conditions \cite{diez2018stochastic}. The integration of UQ methods within stochastic design optimization procedures for vehicle problems was addressed in the AVT-252 group on ``Stochastic Design Optimization for Naval and Aero Military Vehicles,'' where the hull form of a naval destroyer was optimized for stochastic ocean conditions \cite{serani2019stochastic}. The group also addressed the application of several UQ methods (including multi-fidelity approaches) to an airfoil benchmark problem and the results were discussed in \cite{quagliarella2019benchmarking}. Finally, the application of multi-fidelity methodologies to the analysis and design of vehicles is addressed by the AVT-331 group on ``Goal-Driven, Multi-Fidelity Approaches for Military Vehicle System-Level Design.'' An overview on the AVT-331 activities on multi-fidelity approaches may be found in \cite{beran2020}.

There is by now a large consensus in the UQ community on the fact that large-scale, industrially relevant UQ analyses can only be performed by leveraging on multi-fidelity methodologies, i.e., methodologies that explore the bulk of the variability of the quantities of interest of the simulation over coarse meshes (or more generally, computationally inexpensive models with e.g. simplified physics), and resort to querying high-fidelity models (e.g., refined meshes or full-physics models) only sparingly, to correct the initial guess produced with the low-fidelity models. Several approaches to this general framework can be conceived, depending on the kind of fidelity models considered and on the strategy used to sample the parameter space (i.e., for what values of the uncertain parameters the different fidelity models should be queried/evaluated).

In this context, the family of multi-level/multi-index methods has
received an increasing attention, due to its effectiveness and solid
mathematical ground. The hierarchy of models considered by these
methods is usually obtained by successive (dyadic) refinements of a
computational mesh. The multi-level/multi-index distinction arises
from the number of discretization {hyper-}parameters that are considered to
{control the overall discretization of the problem, i.e., how many
  discretization hyper-parameters are used to determine the computational meshes
  (e.g. one or multiple size parameters $h_1$, $h_2$, $h_3$
  for the mesh elements and/or time-stepping) and the number of samples
  from the parameters space to be solved on each mesh (e.g. specified by
  a single number or by a tuple of different numbers along different directions in the parametric space)}.
Combining {these generic discretization strategies} with a
specific sampling strategy over the parameter space results in the
different variations of the method, such as Multi-Level Monte-Carlo
(the first one to be proposed, \cite{giles:MLMC,scheichl.giles:MLMC}),
Multi-Index Monte Carlo \cite{hajiali.eal:MultiIndexMC},
Multi-Level/Multi-Index Quasi-Monte-Carlo \cite{kss12},
Multi-Level Stochastic Collocation \cite{teckentrup.etal:MLSC},
Multi-Index Stochastic Collocation \cite{beck.eal:MISC-IGA,jakeman2019adaptive,hajiali.eal:MISC1,hajiali.eal:MISC2}
Multi-Level Least-Squares polynomial approximation \cite{hajiali2017multilevel}, etc.
{The wording ``Stochastic Collocation'' is to be understood as a synonim of
``sampling in the parametric space'': it refers to the fact that the parameters of the problem can be
seen as random (stochastic) variables, and sampling the parameteric space can be seen as ``collocating the
approximation problem at points of the stochastic domain''.}
The multi-level/multi-index framework can also be extended to
the more generic scenario of Monte-Carlo sampling of multiple
fidelities (e.g. combining different physical models),
see e.g. \cite{peherstorfer:MFsurvey,Gorodetsky2020b}.
See also e.g. \cite{Pisaroni:MLMC-opt,pisaroni:CMLMC,geraci2019} for
applications of Multi-Level Monte-Carlo approaches in the context of
aerodynamics.

Another widely studied class of multi-fidelity methods employs  
%multi-fidelity methods based on 
kernel-based surrogates such as hierarchical kriging \cite{han2012-AIAA}, co-kriging \cite{debaar2015-CF}, Gaussian process \cite{wackers2020-AIAA}, and radial-basis functions \cite{serani2019-IJCFD}. Additive, multiplicative, or hybrid correction methods, also known as ``bridge functions'' or ``scaling functions'' \cite{han2013-AST}, are used to build multi-fidelity surrogates. Further efficiency of multi-fidelity surrogates is gained using dynamic/adaptive sampling strategies, for which the multi-fidelity design of experiments for the surrogate training is not defined a priori but dynamically updated, exploiting the information that becomes available during the training process. Training points are dynamically added with automatic selection of {both} their location and the desired fidelity level, {with} the aim of reducing the computational cost required to properly represent the function. An example of adaptive multi-fidelity  sampling based on the maximum prediction uncertainty is given in \cite{serani2019-IJCFD}.

The objective of the present work is to assess and compare the use of two methods from these two methodological families for the forward UQ of complex industrial problems. Specifically, the performances of Multi-Index Stochastic Collocation (MISC, \cite{beck.eal:MISC-IGA,jakeman2019adaptive,hajiali.eal:MISC1,hajiali.eal:MISC2}) and adaptive Multi-Fidelity Stochastic Radial Basis Functions (SRBF \cite{wackers2020-AIAA}) are compared on the UQ of a roll-on/roll-off passengers (RoPax) ferry sailing in calm water with two operational uncertainties,
specifically ship speed and draught, 
the latter being directly linked to the payload.
The estimation of expected value, standard deviation, and probability density function of the (model-scale) resistance {is} presented and discussed.
{
  Both MISC and SRBF of course need to repeatedly solve the free-surface Navier-Stokes
  equations (i.e. to perform a CFD simulation) for different configurations of the operational parameters.
  The solutions are obtained by the unsteady Reynolds Averaged Navier-Stokes (RANS) solver
  $\chi$navis \cite{dimascio2007,dimascio2009,broglia2018}, developed at CNR-INM. More specifically, both methods
}
 use as fidelity levels the intermediate grids employed by the RANS solver (which is a multi-grid solver):
these grids -- four in total -- are obtained as isotropic {derefinement of an initial fine grid}.
Therefore, both MISC and SRBF are used as multi-index methods with only one component controlling the spatial discretization.

\section{Forward Uncertainty Quantification Method}
% Let us consider an hexahedral mesh with non-cubic elements
% all elements have the same size\footnote{this assumption can be relaxed, but we keep it for simplicity of exposition} and are allowed to be non-cubic, i.e., their edges have size $h_1 = c_1 2^{-\alpha_1}$, $h_2 = c_2 2^{-\alpha_2}$, $h_3 = c_3 2^{-\alpha_3}$, for some constants $c_1, c_2, c_3$ and user-defined integer values $\alpha_1, \alpha_2,\alpha_3$.

Let us consider a single-patch mesh of the computational domain {of a CFD simulation} with non-cubic hexaxedral elements
of the same size\footnote{{The assumption that all elements must be of the same size} can be relaxed, but it is kept for simplicity of exposition} 
and let us also assume that the level of refinement of the mesh along each physical direction can be specified by prescribing 
some user-defined integer values $\alpha_1, \alpha_2,\alpha_3$; to fix ideas, one can think e.g. that the size of each element of the mesh
scales as $2^{-\alpha_1} \times 2^{-\alpha_2} \times 2^{-\alpha_3}$, but this is not necessary. 
The three values of $\alpha_i$ are collected in a multi-index $\bm{\alpha} = [\alpha_1, \alpha_2, \alpha_3]$; prescribing the multi-index $\bm{\alpha}$ thus prescribes the computational mesh to be generated. If this flexibility is not allowed by the mesh-generator (or by the problem itself), it is possible to set $\alpha_1 = \alpha_2 = \alpha_3=\alpha$, i.e., controlling the mesh-generation by a single integer value $\alpha$ {(this is actually the case for the RoPax ferry example considered in this work)}.
The same philosophy applies also to multi-patch meshes, where in principle there could be up to three values $\alpha_i$ for each patch.
The quantity of interest of the CFD simulation computed over the mesh specified by $\bm{\alpha}$ is denoted by $G_{\bm{\alpha}}$; this could be either the full velocity field or a scalar quantity associated to it.

Next, let us assume that the CFD simulation depends on the value of one or more random/uncertain parameters, say $N$ parameters collected in the random vector ${\bf y}=[y_1,y_2,\ldots,y_N]$. Denote by $\Gamma$ the set of all possible values of ${\bf y}$, and by $\rho({\bf y})$ the probability density function (PDF) of the random vector ${\bf y}$ over $\Gamma$. Thus, the primary goal of the forward UQ analysis is to compute an approximation of $\mathbb{E}[G_{\bm{\alpha}}]$, i.e., of the expected value of $G_{\bm{\alpha}}$. This quantity is typically computed by a sampling approach, i.e., the partial differential equation (PDE) at hand is solved (i.e. the CFD simulation is performed) over the grid indexed by $\bm{\alpha}$ for several possible values of ${\bf y}$, and the results are averaged with some weights:
\begin{equation}\label{eq:MC}
\mathbb{E}[G_{\bm{\alpha}}] \approx \sum_{j=1}^J G_{\bm{\alpha}}({\bf y}_j) \omega_j.
\end{equation}

The simplest averaging scheme is Monte Carlo, where the values ${\bf y}_j$ are chosen at random over $\Gamma$ (according to the PDF $\rho$) and $\omega_j = 1/J$. 
Other quantities can be object of the forward UQ analysis, e.g. the standard deviation $\text{std}[G_{\bm{\alpha}}]  =\mathbb{E}[G_{\bm{\alpha}}^2] - \mathbb{E}[G_{\bm{\alpha}}]^2$ or the PDF of $G_{\bm{\alpha}}$. More details on the computation of the latter will be given in Sect. \ref{sect:numerical_results}.

\subsection{Multi-Index Stochastic Collocation (MISC)} \label{sect:misc}
{In this section,  we introduce the MISC method for forward UQ.} 
Roughly speaking, MISC is based on using as {quadrature} points ${\bf y}_j$ {in Eq.\ \eqref{eq:MC}} the union of several Cartesian grids over the domain $\Gamma$, that are obtained by tensorization of univariate quadrature rules (which should be chosen according to $\rho({\bf y})$ for computational efficiency).
In {the RoPax ferry example considered in} this work, $y_1,y_2,\ldots,y_N$ are uniform and independent random variables (see Section \ref{sect:problem_description}),
{therefore} the univariate Clenshaw--Curtis (CC) quadrature is employed. The points for the $K$-points {univariate}
CC quadrature rule are
\begin{equation*}
  t_{j,K}=\cos\left(\frac{(j-1) \pi}{K-1}\right), \quad 1\leq j \leq K,
\end{equation*}
and the corresponding quadrature weights can be efficiently computed by fast Fourier transform, see e.g. \cite{trefethen:comparison}.
Similarly to what done with the multi-index $\bm{\alpha}$ for the physical domain, a multi-index
$\bm{\beta} \in \mathbb{N}^N$ is introduced, that specifies how many points ${\bf y}$ will be used to generate the above-mentioned Cartesian grids.
More specifically, after having introduced the auxiliary function 
\begin{equation}\label{eq:lev_fun}
	m(0)=0, \ m(1)=1, \ m(i) = 2^{i-1}+1 \ \text{for} \ i \geq 2, 
\end{equation} 
$m(\beta_1)$ CC points are generated for $y_1$, $m(\beta_2)$ CC points are generated for $y_2$ etc., and the grid obtained by taking the Cartesian product of the $N$ sets of points thus generated is considered.
Note that this choice of $m$ guarantees that,
  given any two multi-indices $\bm{\beta}_1$ and $\bm{\beta}_2$,
  the grid obtained using $\bm{\beta}_1$ is contained (\textit{nested})
  in the one obtained using $\bm{\beta}_2$ whenever all components of $\bm{\beta}_1$ are smaller or equal than the
  corresponding components in $\bm{\beta}_2$. This is clearly useful in the context of adaptive schemes,
  like the version of MISC that we advocate in this work.
The quadrature weight $\omega_j$ of each point of the Cartesian grid is immediately obtained by 
taking the product of the corresponding univariate weights.

The approximation of $\mathbb{E}[G_{\bm{\alpha}}]$ computed over this grid with Eq. \eqref{eq:MC} is denoted as $\mathcal{Q}_{\bm{\alpha},\bm{\beta}}$. Clearly, it would be beneficial to have both multi-indices $\bm{\alpha}$ and $\bm{\beta}$ with large components, say $\bm{\alpha} = \bm{\alpha}^\star$ and $\bm{\beta} = \bm{\beta}^\star$, i.e., to average the values of many PDE solutions over a refined computational mesh. However, this is typically unfeasible due to computational costs. One possible remedy is to exploit the fact that a single, highly refined approximation $\mathcal{Q}_{\bm{\alpha}^\star,\bm{\beta}^\star}$ can often%
\footnote{whenever $G({\bf y})$ is a smooth function with respect to ${\bf y}$, i.e., {roughly speaking, whenever}
  small changes in ${\bf y}$ imply small changes in $G({\bf y})$}
be approximated as a linear combination of many coarser $\mathcal{Q}_{\bm{\alpha},\bm{\beta}}$, where whenever the spatial discretization $\bm{\alpha}$ is refined, the quadrature level $\bm{\beta}$ is kept to a minimum and viceversa (of course, the combined cost of computing the set of coarse discretizations should be smaller than the cost of the highly refined one). This is in a nutshell the idea of MISC. In other words, MISC is a classical multi-level scheme, where most of the statistical variability of $G$ is explored by solving many PDEs with coarse meshes (large $\|\bm{\beta}\|$ with small $\|\bm{\alpha}\|$)\footnote{$\|\cdot\|$ denotes the Euclidean norm}
and then the result is corrected with a few PDE solutions with refined meshes (large $\|\bm{\alpha}\|$ with small $\|\bm{\beta}\|$).
Clearly, as all multi-level approaches, MISC works well only if the levels are sufficiently separated,
i.e. if the number of degrees of freedom {of the computational mesh}
(and the corresponding computational cost) grows significantly from one level to the next one:
to fix ideas again, things will work {well}
if the number of elements in the mesh scales e.g. as $2^{\alpha_1} \times 2^{\alpha_2} \times 2^{\alpha_3}$,
but not if e.g. increasing $\alpha_1$ to $\alpha_1+1$ adds only one element to the mesh.

\begin{figure}[!t]
  \centering
  \subfigure[]{\includegraphics[width=0.32\linewidth]{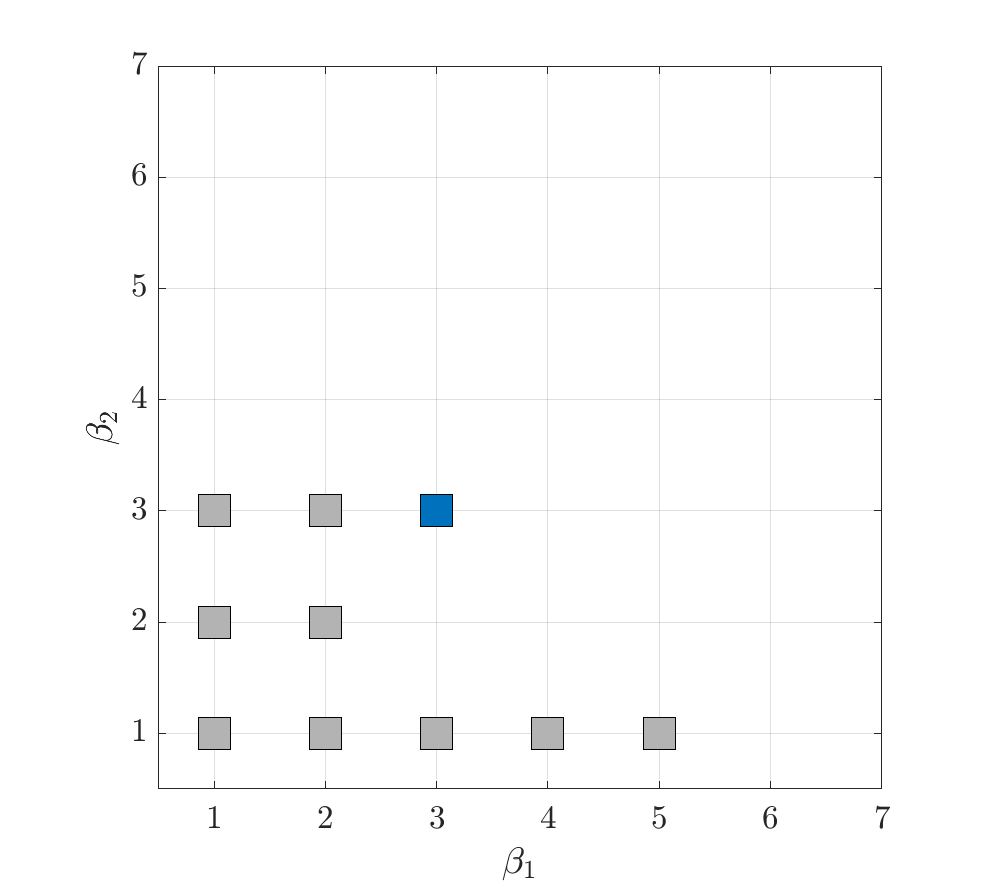}}
  \subfigure[]{\includegraphics[width=0.32\linewidth]{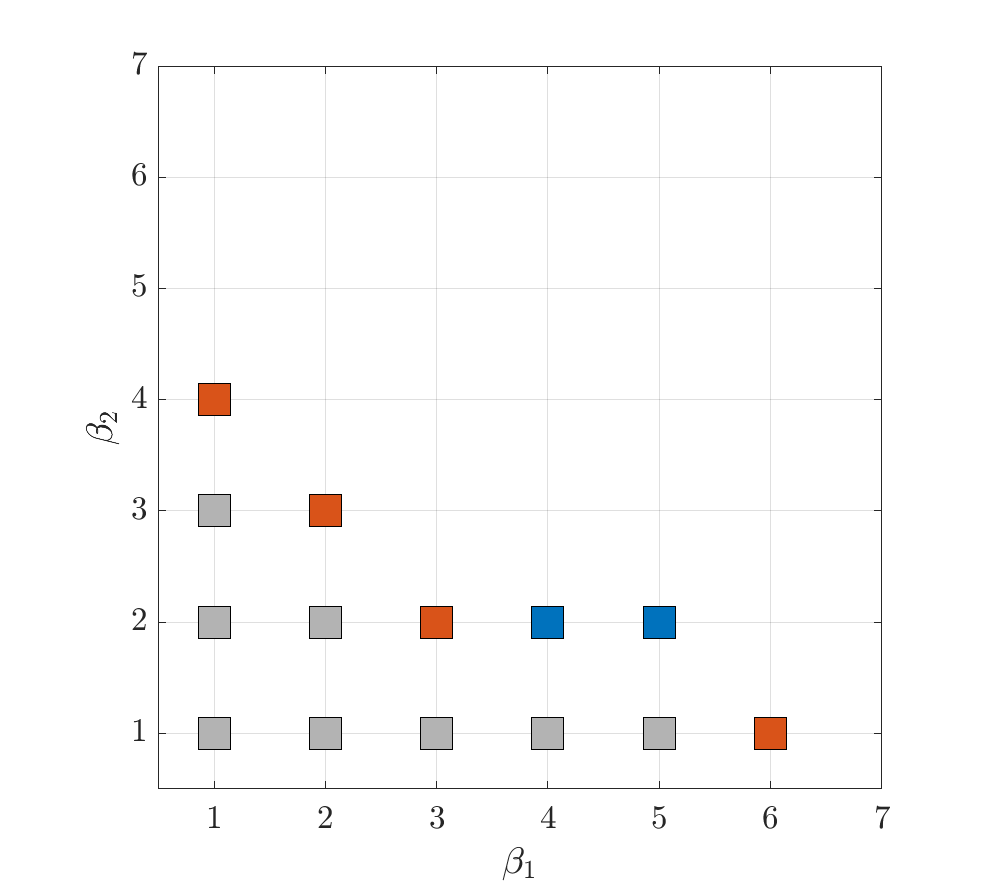}}
  \caption[index_sets]{(a): the gray set is downward closed, whereas adding the blue multi-index to it would result
    in a set not downward closed; (b): a downward closed set (in gray) and its margin (indices marked in red and blue).
    If Algorithm \ref{algo:misc_implementation} reaches the gray set, it will next explore all indices marked in red (their addition to the gray
    set keeps the downward closedness property) but not those marked in blue. The red set is also known as ``reduced margin''. }
  \label{fig:index_sets}
\end{figure}

In the following the details of MISC are briefly recalled, following closely the setup in \cite{beck.eal:MISC-IGA}. Firstly, the so-called univariate and multivariate ``detail operators'' on the physical and {parametric} domains have to be introduced. They are defined as follows, with the understanding that $\mathcal{Q}_{\bm{\alpha},\bm{\beta}}=0$ when at least one component of $\bm{\alpha}$ or $\bm{\beta}$ is zero.
{Also, for notational consistence with the existing MISC literature,
  in the rest of this section and in the formulas we drop the word ``parametric'' and use the word ``stochastic'' instead,
  with the understanding that they are synonims to our purposes:}
\begin{alignat*}{2}
&\text{\textbf{Univariate physical detail: }}
&& \Delta_i^{\text{phys}}[\mathcal{Q}_{\bm{\alpha},\bm{\beta}}]=\mathcal{Q}_{\bm{\alpha},\bm{\beta}}-\mathcal{Q}_{\bm{\alpha}-\ee_i,\bm{\beta}} \text{ with } 1 \leq i \leq 3; \\
&\text{\textbf{Univariate stochastic detail: }}
&& \Delta_i^{\text{stoc}}[\mathcal{Q}_{\bm{\alpha},\bm{\beta}}]=\mathcal{Q}_{\bm{\alpha},\bm{\beta}}-\mathcal{Q}_{\bm{\alpha},\bm{\beta}-\ee_i} \text{ with } 1 \leq i \leq N; \\
&\text{\textbf{Multivariate physical detail: }}
&& \bm{\Delta}^{\text{phys}}[\mathcal{Q}_{\bm{\alpha},\bm{\beta}}] = \bigotimes_{i=1}^3 \Delta_i^{\text{phys}}[\mathcal{Q}_{\bm{\alpha},\bm{\beta}}]; \\% with $1 \leq i \leq 3$;
&\text{\textbf{Multivariate stochastic detail: }}
&& \bm{\Delta}^{\text{stoc}}[\mathcal{Q}_{\bm{\alpha},\bm{\beta}}] = \bigotimes_{j=1}^N \Delta_j^{\text{stoc}}[\mathcal{Q}_{\bm{\alpha},\bm{\beta}}]; \\% with $1 \leq i \leq N$;
&\text{\textbf{Mixed multivariate detail: }}
&& \bm{\Delta}^{\text{mix}}[\mathcal{Q}_{\bm{\alpha},\bm{\beta}}] = \bm{\Delta}^{\text{stoc}}\left[\bm{\Delta}^{\text{phys}}[\mathcal{Q}_{\bm{\alpha},\bm{\beta}}] \right].  
\end{alignat*}
Note that $\ee_i$ denotes the canonical multi-index, i.e. $(\ee_i)_k = 1$ if $i=k$ and 0 otherwise.
Observe that taking tensor products of univariate details amounts to composing their actions, e.g.,
\[
\bm{\Delta}^{\text{phys}}[\mathcal{Q}_{\bm{\alpha},\bm{\beta}}]
= \bigotimes_{i=1}^3 \Delta_i^{\text{phys}}[\mathcal{Q}_{\bm{\alpha},\bm{\beta}}]
= \Delta_1^{\text{phys}}\left[ \, \Delta_2^{\text{phys}}\left[\Delta_3^{\text{phys}}\left[ \mathcal{Q}_{\bm{\alpha},\bm{\beta}} \right] \, \right] \, \right], %\bigg]
\]
and analogously for the stochastic multivariate detail operators, $\bm{\Delta}^{\text{stoc}}[\mathcal{Q}_{\bm{\alpha},\bm{\beta}}]$.
Crucially, this in turn implies that the multivariate operators can be evaluated by evaluating certain full-tensor approximations
$\mathcal{Q}_{\bm{\alpha},\bm{\beta}}$ and then taking linear combinations:
\begin{align*}
\bm{\Delta}^{\text{phys}}[\mathcal{Q}_{\bm{\alpha},\bm{\beta}}]
& = \Delta_1^{\text{phys}}\left[ \, \Delta_2^{\text{phys}}\left[ \Delta_3^{\text{phys}}\left[ \mathcal{Q}_{\bm{\alpha},\bm{\beta}} \right] \, \right] \, \right] %\bigg]
= \sum_{\jj \in \{0,1\}^3} (-1)^{|\jj|} \mathcal{Q}_{\bm{\alpha}-\jj,\bm{\beta}};\\
\bm{\Delta}^{\text{stoc}}[\mathcal{Q}_{\bm{\alpha},\bm{\beta}}]
&  = \sum_{\jj \in \{0,1\}^N} (-1)^{|\jj|} \mathcal{Q}_{\bm{\alpha},\bm{\beta}-\jj}.
\end{align*}
The latter expressions are known in the sparse-grids community as ``combination-technique'' formulations, and can be very
useful for practical implementations. In particular, they allow to evaluate e.g. $\bm{\Delta}^{\text{phys}}[\mathcal{Q}_{\bm{\alpha},\bm{\beta}}]$
by calling pre-existing softwares on different meshes up to $2^3$ times in a ``black-box'' fashion.
Analogously, evaluating  $\bm{\Delta}^{\text{stoc}}[\mathcal{Q}_{\bm{\alpha},\bm{\beta}}]$
requires evaluating up to $2^N$ operators $\mathcal{Q}_{\bm{\alpha},\bm{\beta}}$ over different quadrature grids, and
evaluating $\bm{\Delta}^{\text{mix}}[\mathcal{Q}_{\bm{\alpha},\bm{\beta}}]$ requires evaluating up to $2^{3+N}$ operators $\mathcal{Q}_{\bm{\alpha},\bm{\beta}}$ over different quadrature grids and physical meshes.
Observe that by introducing these detail operators a hierarchical decomposition of $\mathcal{Q}_{\bm{\alpha},\bm{\beta}}$ is available; indeed, the following is a telescopic identity
\begin{equation}\label{eq:telescopic_sum}
\mathcal{Q}_{\bm{\alpha},\bm{\beta}} = \sum_{[\bm{\ii},\bm{\jj}] \leq [\bm{\alpha},\bm{\beta}]} \bm{\Delta}^{\text{mix}}[\mathcal{Q}_{\bm{\ii},\bm{\jj}}],
\end{equation}
i.e., it can be easily verified by replacing each term $\bm{\Delta}^{\text{mix}}[\mathcal{Q}_{\bm{\ii},\bm{\jj}}]$ with the
the corresponding combination-technique formula that all terms except $\mathcal{Q}_{\bm{\alpha},\bm{\beta}}$ will cancel.
As an example, the case of one-dimensional physical and stochastic spaces can be considered. Recalling that by definition $\mathcal{Q}_{i,j} = 0$ when either $i=0$ or $j=0$, it can be seen that
\begin{align*}
\sum_{[i,j] \leq [2,2]} \bm{\Delta}^{\text{mix}}[\mathcal{Q}_{i,j}] 
& = \bm{\Delta}^{\text{mix}}[\mathcal{Q}_{1,1}]
+ \bm{\Delta}^{\text{mix}}[\mathcal{Q}_{1,2}]
+ \bm{\Delta}^{\text{mix}}[\mathcal{Q}_{2,1}]
+ \bm{\Delta}^{\text{mix}}[\mathcal{Q}_{2,2}] \\
& = \mathcal{Q}_{1,1}
+ ( \mathcal{Q}_{1,2} - \mathcal{Q}_{1,1} )
+ ( \mathcal{Q}_{2,1} - \mathcal{Q}_{1,1} )
+ \mathcal{Q}_{2,2} - \mathcal{Q}_{2,1} - \mathcal{Q}_{1,2} - \mathcal{Q}_{1,1} \nonumber \\
& = \mathcal{Q}_{2,2}. \nonumber
\end{align*}
The crucial observation is that not all of the details in the above hierarchical decomposition \eqref{eq:telescopic_sum} contribute equally to the approximation,
i.e., {some of them} can be discarded and the resulting formula will retain good approximation properties
at a fraction of the computational cost. Thus, the MISC approximation of $\mathbb{E}[G_{\bm{\alpha}}]$ is introduced as
\[
\mathbb{E}[G_{\bm{\alpha}}] \approx \mathcal{Q}_{\Lambda}=\sum_{[\bm{\alpha},\bm{\beta}] \in \Lambda} \bm{\Delta}^{\text{mix}}[\mathcal{Q}_{\bm{\alpha},\bm{\beta}}],
\]
for a suitable multi-index set $\Lambda \subset \mathbb{N}^{3+N}$,
which should be chosen as downward closed, i.e.  (see also Fig. \ref{fig:index_sets}a)
\[
\forall \mathbf{k} \in \Lambda, \quad \mathbf{k} - \ee_j \in \Lambda \text{ for every } j=1,\ldots,3+N \text{ such that } \mathbf{k}_j > 1.  
\]
Clearly, the MISC estimator has a combination-technique expression as well, which can be written in compact form as
\begin{equation}\label{eq:misc}
\mathbb{E}[G_{\bm{\alpha}}] \approx  \mathcal{Q}_\Lambda
=\sum_{[\bm{\alpha},\bm{\beta}] \in \Lambda} \bm{\Delta}^{\text{mix}}[\mathcal{Q}_{\bm{\alpha},\bm{\beta}}]
=\sum_{[\bm{\alpha},\bm{\beta}] \in \Lambda} \sum_{\substack{[\bm{i},\bm{j}] \in \{0,1\}^{3+N}\\ [\bm{\alpha}+\bm{i},\bm{\beta}+\bm{j}] \in \Lambda}}
(-1)^{\lvert [\bm{i},\bm{j}] \rvert} \mathcal{Q}_{\bm{\alpha},\bm{\beta}},
\end{equation}
which finally shows the initial statement that the MISC evaluation is computed by evaluating full-tensor operators $\mathcal{Q}_{\bm{\alpha},\bm{\beta}}$ independently 
and combining them linearly according to \eqref{eq:misc}. This is the approximation formula which is used in {our} pratical implementation of the MISC method.

Of course, the effectiveness of the MISC estimator depends on the choice of the multi-index set $\Lambda$.
Several strategies have been explored in the literature; a suitable set $\Lambda$ can either be designed a-priori, by a careful analysis of the PDE at hand, see e.g. \cite{beck.eal:MISC-IGA}, or on-the-run by adaptive algorithms, see e.g. \cite{jakeman2019adaptive}; in this contribution the focus is on the latter option.
To devise an optimal strategy for selecting a good multi-index set, the following error decomposition is introduced
\begin{align*}
\lvert \mathbb{E}[G_{\bm{\alpha}}] - \mathcal{Q}_\Lambda \rvert
&  = \Big \lvert \mathbb{E}[G_{\bm{\alpha}}] - \sum_{[\bm{\alpha},\bm{\beta}] \in \Lambda} \bm{\Delta}^{\text{mix}}[\mathcal{Q}_{\bm{\alpha},\bm{\beta}}] \Big\rvert \nonumber \\
&  = \Big \lvert \sum_{[\bm{\alpha},\bm{\beta}] \not \in \Lambda} \bm{\Delta}^{\text{mix}}[\mathcal{Q}_{\bm{\alpha},\bm{\beta}}] \Big \rvert
\leq \sum_{[\bm{\alpha},\bm{\beta}] \not \in \Lambda} \big \lvert \bm{\Delta}^{\text{mix}}[\mathcal{Q}_{\bm{\alpha},\bm{\beta}}] \big \rvert
= \sum_{[\bm{\alpha},\bm{\beta}] \not \in \Lambda} E_{\bm{\alpha},\bm{\beta}},  
\end{align*}
where $E_{\bm{\alpha},\bm{\beta}} := \big \lvert \bm{\Delta}^{\text{mix}}[\mathcal{Q}_{\bm{\alpha},\bm{\beta}}] \big \rvert$;
$E_{\bm{\alpha},\bm{\beta}}$ thus represents the ``error contribution'' of $[\bm{\alpha},\bm{\beta}]$,
i.e., the reduction in the approximation error due to having added $[\bm{\alpha},\bm{\beta}]$ to the current
index-set $\Lambda$; in formula 
\begin{equation}\label{eq:error_contr}
E_{\bm{\alpha},\bm{\beta}} = \lvert \mathcal{Q}_{\Lambda \cup [\bm{\alpha},\bm{\beta}]} -\mathcal{Q}_{\Lambda} \rvert.
\end{equation}
Similarly, the ``work contribution'' $W_{\bm{\alpha},\bm{\beta}}$ is defined as the work required
to add $[\bm{\alpha},\bm{\beta}]$ to the current index-set $\Lambda$. It is the product of the computational cost associated to the spatial grid identified by the multi-index $\bm{\alpha}$, denoted by $\text{cost}(\bm{\alpha})$ (see details in Sect. \ref{sect:numerical_results}, Eq. \eqref{eq:cost}), times the number of new evaluations of the PDE required by the multi-index $\bm{\beta}$, i.e.  
\begin{equation}\label{eq:work_contr}
W_{\bm{\alpha},\bm{\beta}} = \text{cost}(\bm{\alpha})  \prod_{n=1}^N (m(\beta_n)-m(\beta_n-1)),
\end{equation}
with $m$ defined as in \eqref{eq:lev_fun}.
Note that the expression above is based on the fact that the employed quadrature rule is nested.  
An effective strategy to build adaptively a MISC approximation can then be broadly described as follows: given the MISC approximation
associated to the multi-index set $\Lambda$, a new MISC approximation is built by adding to $\Lambda$ the multi-index with the next
largest profit $P_{\bm{\alpha},\bm{\beta}} = \frac{E_{\bm{\alpha},\bm{\beta}}}{W_{\bm{\alpha},\bm{\beta}}}$.
Of course, it is impossible to compute in advance the profits of all multi-indices in $\mathbb{N}^{3+N}$,
and {we also need} to guarantee that $\Lambda$ is downward-closed
at all times. Therefore, in practice the implementation reported in Algorithm \ref{algo:misc_implementation} is used.
It makes use of an auxiliary multi-index set, i.e. the margin of a multi-index set $\Lambda$, $\text{Mar}(\Lambda)$,
which is defined as the set of multi-indices that can be reached ``within one step'' from $\Lambda$ (see also Fig. \ref{fig:index_sets}b)
\begin{equation*}
\text{Mar}(\Lambda) = \{\ii \in \mathbb{N}^{3+N} \text{ s.t. } \ii = \jj + \ee_k \text{ for some } \jj \in \Lambda \text{ and some } k \in \{1,\ldots,3+N\} \}.    
\end{equation*}
This algorithm was first proposed in the context of sparse-grids quadrature in \cite{gerstner.griebel:adaptive} and its MISC implementation was first proposed in \cite{jakeman2019adaptive}.

Note that the MISC method can also be used to obtain a surrogate model for $G_{\bm{\alpha}}({\bf y})$. More specifically, a formula analogous to \eqref{eq:misc} can be written replacing the tensor quadrature operator $\mathcal{Q}_{\bm{\alpha}, \bm{\beta}}$ by a tensor interpolation operator using global Lagrange polynomials collocated at the CC points, see  e.g. \cite{jakeman2019adaptive}. The result is a linear combination of tensor interpolants according to the combination technique formula
{\eqref{eq:misc}}, which is used as surrogate response surface in Sect. \ref{sect:numerical_results}.
A further remark is that the choice of computing the profit based on the ``improvement'' of the approximation of the expected value of $G$ is somehow arbitrary, and many alternatives can be considered, including interpolation-based versions, see e.g.
\cite{nobile.eal:adaptive-lognormal,klimke:thesis,chkifa:adaptive-interp} for the analogous discussion in the context of sparse-grids quadrature
and interpolation. 

Similarly, many stopping criteria can be considered (and possibly used simultaneously), which typically check that computational work,
error contributions or profit estimator are below a desired threshold. More sophisticated error estimators are subject of research,
see e.g. \cite{Guignard:a-post} in the context of sparse grids quadrature.

\begin{algorithm}
	\Fn{\AlCapSty{Multi-Index Stochastic Collocation}}{
		$\bm{I} = \{ (\bm{1},\bm{1})\}$, $\bm{G} = \{(\bm{1},\bm{1})\}$, $\bm{R}_{\bm{I}} = \emptyset$ \tcp*[f]{{Here $\bm{G}$ is a set, not the quantity of interest of the CFD simulation; we use $\bm{G}$ nonetheless, for consistence with previous lit.}}\; 
		Compute MISC estimate $\mathcal{Q}_{G}$ as in \eqref{eq:misc} \;
		\While{ stopping criteria are not met }{
                  \For(\tcp*[f]{for short, $\jj = [\bm{\alpha} , \bm{\beta}]$.}){ %See also Figure \ref{fig:index_sets}
                    $\jj \in \text{Mar}(\bm{I})$ {\bf and} $\bm{I} \cup \{\jj\}$ downward closed }{ 
				Compute MISC estimate $\mathcal{Q}_{ \bm{G}\cup \{\jj\} }$ as in \eqref{eq:misc} \;
				Compute error contribution $E_{\jj}$ as in \eqref{eq:error_contr} \;
				Compute work contribution $W_{\jj}$ as in \eqref{eq:work_contr} \;
				Compute profit $P_{\jj} = E_{\jj} / W_{\jj}$ \;
				$ \bm{G} = \bm{G}\cup \{\jj\}$, $\bm{R}_{\bm{I}} = \bm{R}_{\bm{I}} \cup \{\jj\}$ \;
		}
	Choose $\ii \in \bm{R}_{\bm{I}}$ with the highest profit \;
	$\bm{I} = \bm{I} \cup \{\ii\}$,  $\bm{R}_{\bm{I}} = \bm{R}_{\bm{I}} \setminus \{\ii\}$ \;
	} % end While
	}% end \Fn
	\caption{MISC implementation}
	\label{algo:misc_implementation}
\end{algorithm}

\subsection{Multi-Fidelity Stochastic Radial Basis Functions (SRBF)}
{As an alternative to MISC, we present in this section a methodology for forward UQ
  based on a multi-fidelity surrogate model built on kernel functions (Stochastic Radial Basis Functions).
  Here ``Stochastic'' denotes not only the fact that we are sampling parameters that are random variables, but also to the fact that
  the SRBF method treats one of its hyper-parameters as a random variable, as will be clearer later on.
  Note that for ease of notation in this section we consider that the CFD mesh generation is
  controlled by a scalar value $\alpha$, i.e., the quantity of interest is denoted by $G_\alpha$, $\alpha = 1,\ldots, M$.}

\begin{figure}[!t]
\centering
\includegraphics[width=0.495\textwidth]{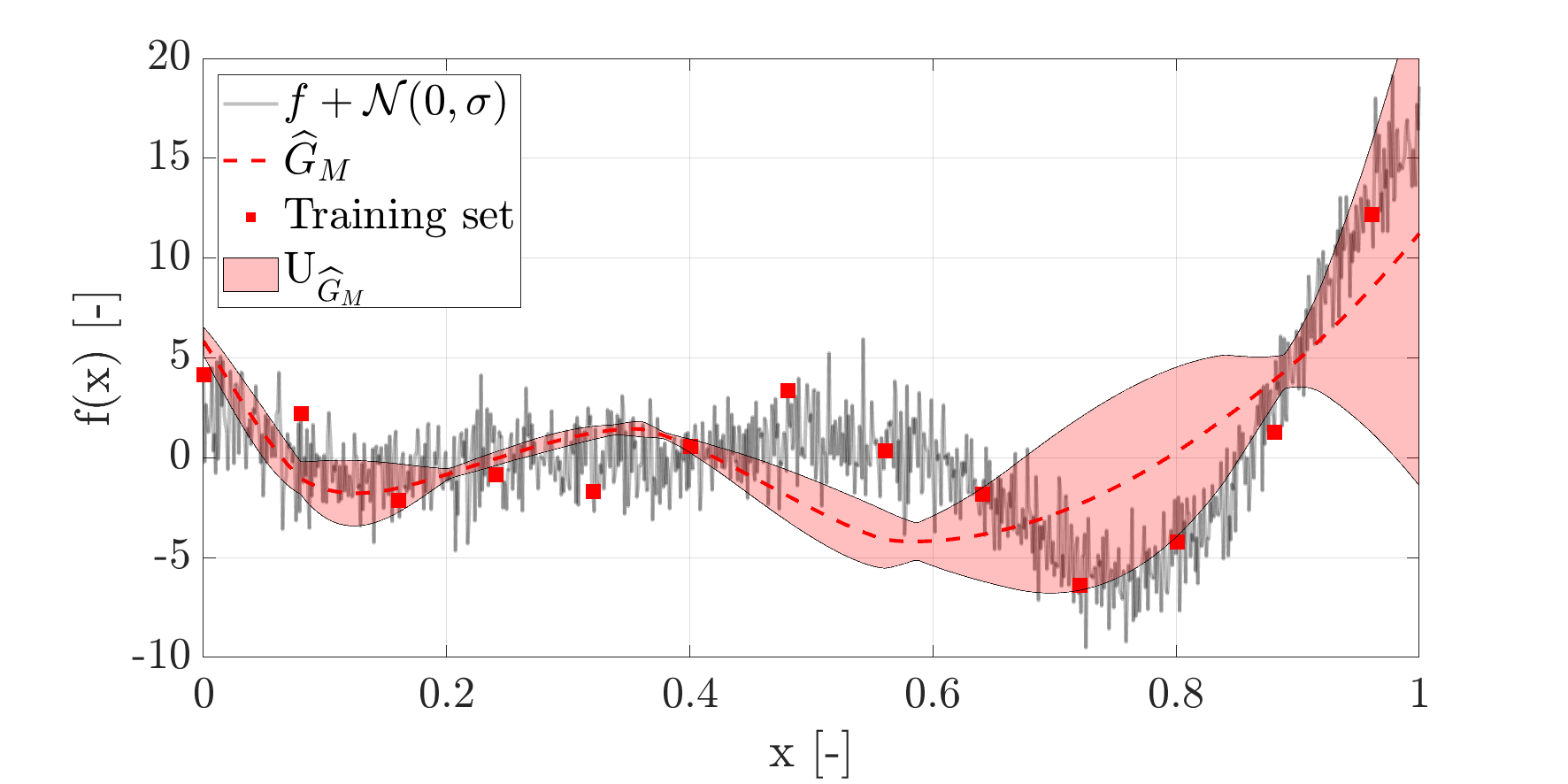}
\caption{SRBF example with least-squares regression.}\label{fig:lsrbf}
\end{figure}

Given a training set {$\mathcal{T}_\alpha=\{ {\bf y}_i,G_\alpha({\bf y}_i) \}_{i=1}^{\mathcal{J}_\alpha}$} 
and normalizing the uncertain parameters domain into a unit hypercube, the RBF prediction is (here) based on a power function kernel and reads
\begin{equation}
f\left({\bf y},\tau\right)=\sum^{\mathcal{K}}_{j=1} w_j ||{\bf y}-{\bf c}_j||^{\tau},
\end{equation}
where $w_j$ are unknown coefficients, ${\bf c}_j$ are {$\mathcal{K}$ points in $\Gamma$ called} RBF centers,
and $\tau\sim \textrm{unif}[\tau_{\min},\tau_{\max}]$ is a stochastic tuning parameter that follows a uniform distribution. RBF models have been widely applied in engineering problems using linear ($\tau=1$, providing a polyharmonic spline of first order \cite{gutmann2001}) and cubic ($\tau=3$, providing a polyharmonic spline of third order \cite{forrester2009}) kernels. This suggests the range of $\tau$ to be defined within $\tau_{\min}=1$ and $\tau_{\max}=3$. Note that the choice of the distribution for $\tau$ is arbitrary and, from a Bayesian viewpoint, this represents the degree of belief in the definition of the tuning parameter.   
The SRBF {surrogate model $\widetilde{G}_\alpha\left({\bf y}\right)$} is computed as the expected value (approximated by {Monte Carlo})
of $f$ over $\tau$ \citep{volpi2015-SMO}:
\begin{eqnarray}\label{eq:SRBF}
{\widetilde{G}_\alpha\left({\bf y}\right)}=\mathbb{E}\left[f\left({\bf y},\tau\right)\right]_\tau\approx\frac{1}{\Theta}\sum\limits_{i=1}^{\Theta}f\left({\bf y},\tau_i\right),
\end{eqnarray}
where $\Theta$ is the number of samples for $\tau$, here {set} equal to $1000$.
{To give more flexibility to the method,}
the coordinates of the RBF centers ${\bf c}_j$ are not {a-priori set to be} coincident with the training points, but {rather chosen by a} $k$-means clustering {algorithm applied to the training points, see} \citep{lloyd1982-IEEE}.
{Several values of $\mathcal{K}$ (i.e., number of centers) are tested and their} optimal number
{ $\mathcal{K}^*_\alpha$ is chosen} by minimizing a leave-one-out cross-validation (LOOCV) metric \citep{li2017-SMO}.
{In details, letting $\tilde{g}_{i,\mathcal{K}}({\bf y})$ be a surrogate model with $\mathcal{K}$ centers
  trained on the whole training set $\mathcal{T}_\alpha$ but the $i$-th point, $\mathcal{K}^*_\alpha$ is defined as:}
\begin{equation}\label{eq:loocv}
{\mathcal{K}^*_\alpha} = {\underset{{\mathcal{K}\in {\mathcal C}}}{\rm argmin}} \,\mathrm{RMSE}({\mathcal{K}}), 
\end{equation}
where $\mathcal C\subset \mathbb{N}$ {and RMSE($\mathcal{K}$) is the root mean square error between each of the
  $\mathcal{J}$ leave-one-out models with $\mathcal{K}$ centers, $\tilde{g}_{1,\mathcal{K}},\ldots,\tilde{g}_{\mathcal{J},\mathcal{K}}$,
  and the evaluations of the full $G_\alpha$ at the point that is being left out for each $\tilde{g}_{i,\mathcal{K}}$:}
\begin{equation}\label{eq:rmse}
  \mathrm{RMSE}({\mathcal{K}}) =
  \sqrt{\dfrac{1}{\mathcal{J}_{{\alpha}}}\sum_{i=1}^{\mathcal{J}_{{\alpha}}}
    \left(G_{{\alpha}}(\mathbf{y}_i) - \tilde{g}_{{i,\mathcal{K}}}(\mathbf{y}_i)\right)^2 },
  \quad {{\bf y}_j \in \mathcal{T}_\alpha.}
\end{equation}
{Clearly, once the optimal number of centers is chosen, the whole set of points is used for the construction of the final surrogate model.}
Whenever the number of RBF centers is lower than the training set size ($\mathcal{K} < {\mathcal{J}_\alpha}$),
the coefficients $w_j$ are determined through a least-squares regression by solving 
\begin{equation}\label{eq:lsrbf}
{\bf w}=\left( {\bf A}^{\mathsf{T}} {\bf A} \right)^{-1} {\bf A}^{\mathsf{T}} {\bf f},
\end{equation}
with ${\bf w}=\{w_j\}$, $a_{ij}=||{\bf y}_i-{\bf c}_j||^{\tau}$ and
{${\bf f}=\{G_\alpha({\bf y}_i)\}$;}
otherwise when the optimal number of RBF centers equals the training set size, exact interpolation at the training points
{($f({\bf y}_i,\tau)=G_\alpha({\bf y}_i)$)}
is imposed and Eq. \eqref{eq:lsrbf} reduces to
\begin{equation}
{\bf Aw}={\bf f},
\end{equation}   
with ${\bf c}_j={\bf y}_j$.
{Having less RBF centers than training points and employing the least-squares approach in Eq. \eqref{eq:lsrbf}
  to determine the coefficients $w_j$ is particularly helpful when the training data are affected by noise.}
An example of least-squares regression is shown in Fig. \ref{fig:lsrbf}.

The uncertainty {$U_{\widetilde{G}_\alpha}\left({\bf y}\right)$}
associated with the SRBF {surrogate model} prediction is
quantified by the 95\%-confidence band of the cumulative density function (CDF) of $f({\bf y}, \tau)$ as follows
\begin{equation}\label{eq:USRBF}
{U_{\widetilde{G}_\alpha}}({\bf y})={\rm CDF}^{-1}(0.975;{\bf y})-{\rm CDF}^{-1}(0.025;{\bf y}),
\end{equation}
with
\begin{equation}
{\rm CDF}(\lambda;{\bf y})=\frac{1}{\Theta}\sum\limits_{i=1}^{\Theta} H[\lambda-f({\bf y},\tau_i)],
\end{equation}
where $H(\cdot)$ is the Heaviside step function.

\begin{figure}[!t]
\centering
\includegraphics[width=0.495\textwidth]{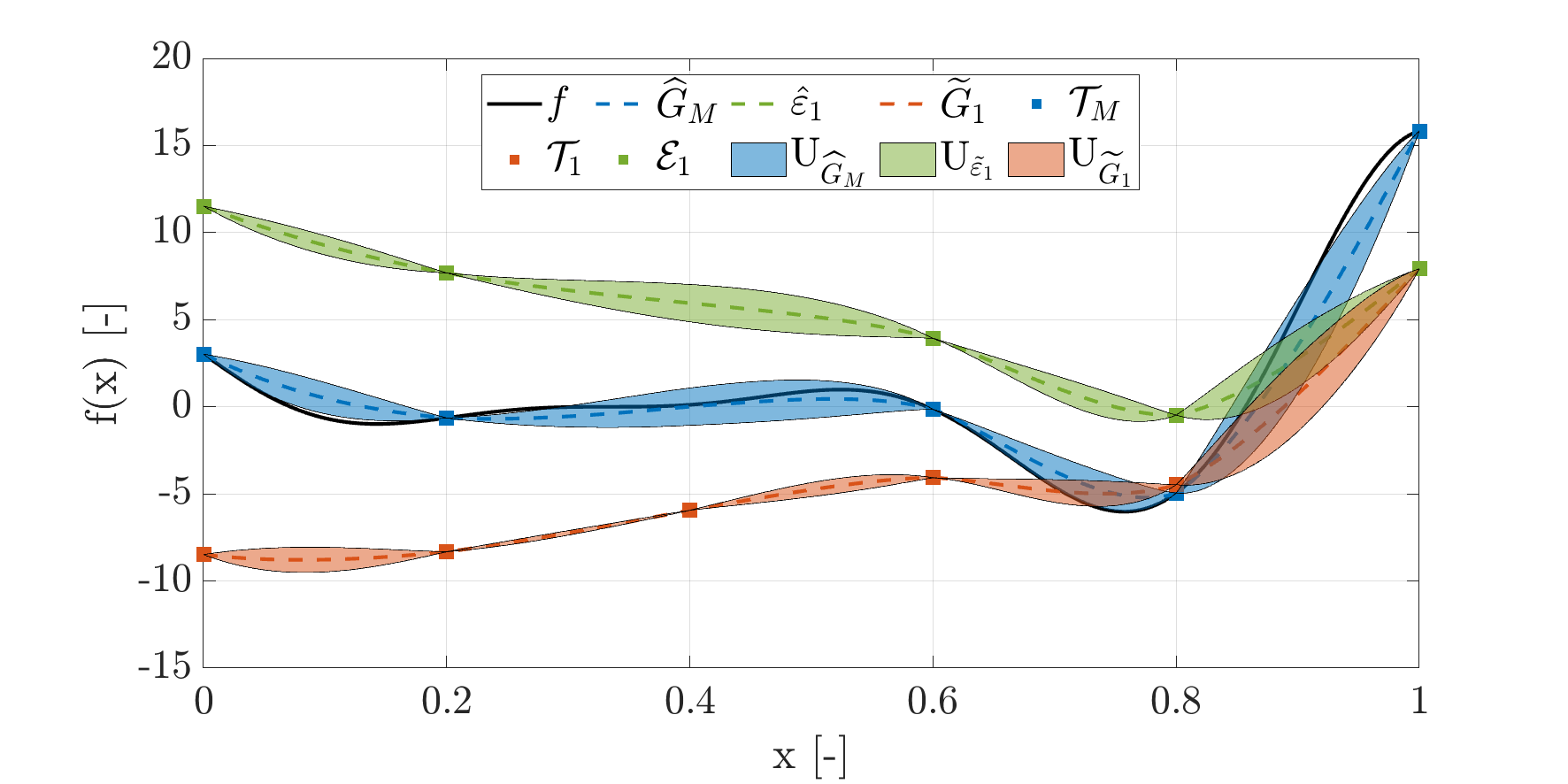}
\caption{Example of multi-fidelity surrogate with $M=2$ and exact interpolation at the training points.}\label{fig:mfm}
\end{figure}

The multi-fidelity approximation {of $G_\alpha$},
introduced in \cite{serani2019-MARINE} and extended to noisy data in \cite{wackers2020-SNH},
is {then} adaptively built as follows.  Extending the definition of the surrogate
{model} training set to an arbitrary number of $M$ fidelity levels as
{$\{\mathcal{T}_{\alpha}\}_{\alpha=1}^M$, with each $\mathcal{T}_{\alpha} = \{{\bf y}_j,G_{\alpha}({\bf y}_j)\}_{j=1}^{\mathcal{J}_{\alpha}}$},
the multi-fidelity approximation $\widehat{G}_{\alpha}(\mathbf{y})$ of $G_{\alpha}(\mathbf{y})$ reads 
\begin{equation}
\label{eq:MLGeneral}
\widehat{G}_{\alpha}(\mathbf{y})\approx\widetilde{G}_1(\mathbf{y})+\sum_{i=1}^{{\alpha}-1}\tilde{\varepsilon}_i(\mathbf{y}), %\\
\end{equation}
where $\tilde{{\varepsilon}}_i(\mathbf{y})$ is the inter-level error surrogate with an associate training set $\mathcal{E}_i=\{({\bf y},G_{i+1}({\bf y})-\widehat{G}_{i}({\bf y}))\,|\,{\bf y} \in \mathcal{T}_{i+1} \cap \mathcal{T}_{i} \}$. It can be noted that Eq. \eqref{eq:MLGeneral} does not strictly require nested training sets. An example of the multi-fidelity approximation with two fidelities is shown in Fig. \ref{fig:mfm}.
\begin{figure}[!b]
\centering
\subfigure[]{\includegraphics[width=0.495\textwidth]{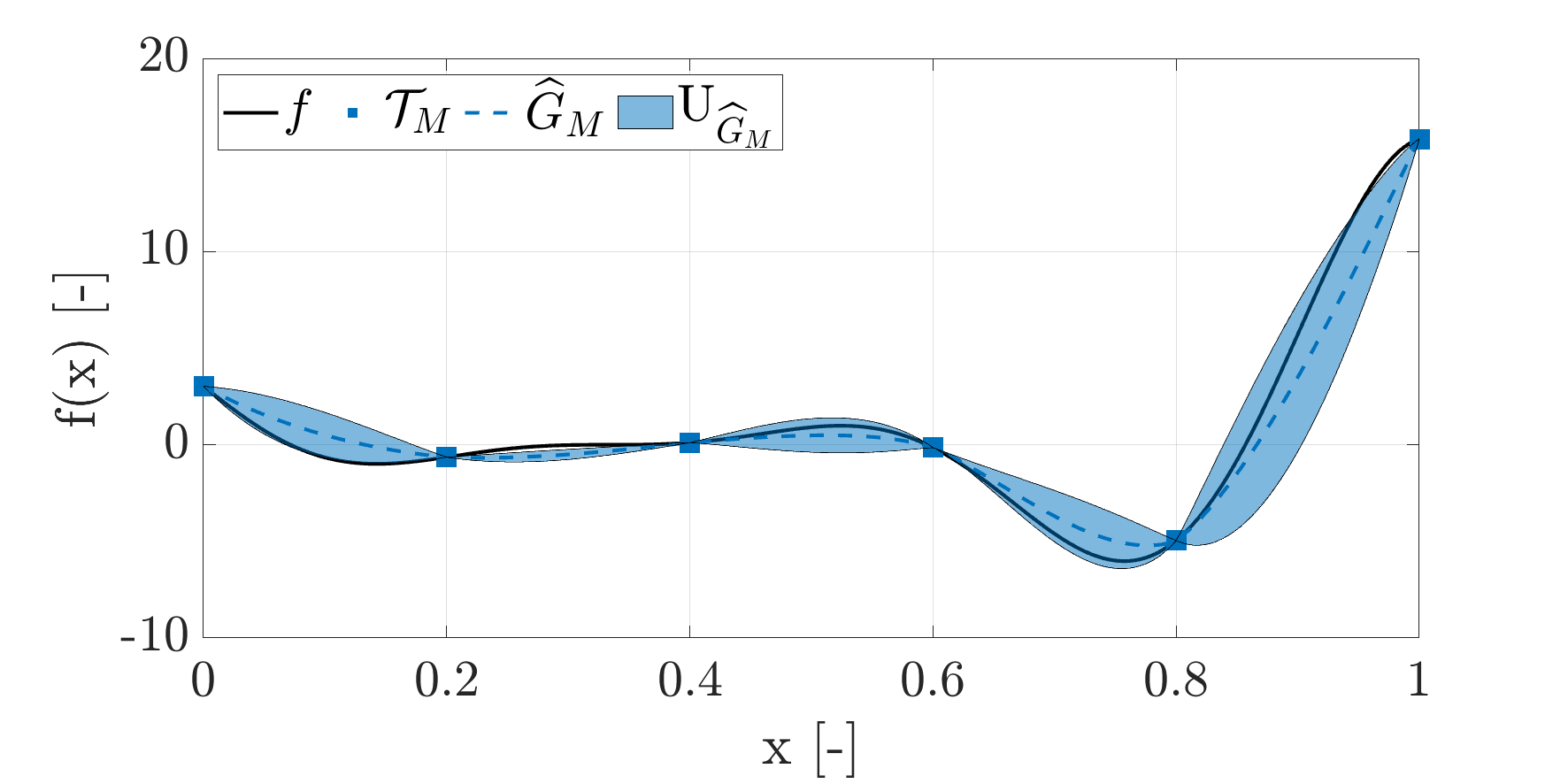}}
\subfigure[]{\includegraphics[width=0.495\textwidth]{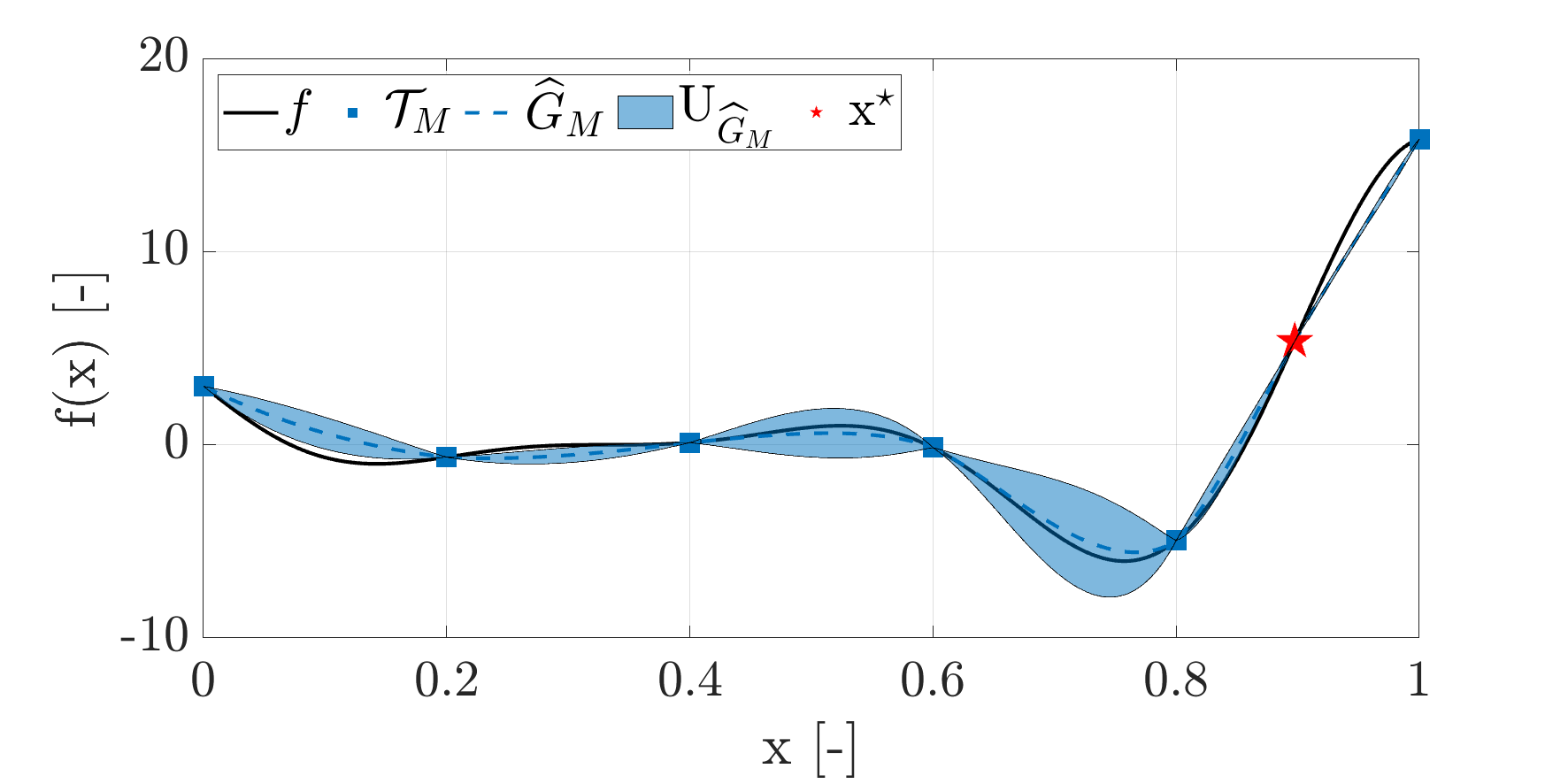}}
\caption{Example of the adaptive sampling method using one fidelity without noise: (a) shows the initial SRBF with the associated prediction uncertainty and training set; (b) shows the position of the new training point, the new SRBF prediction, and its uncertainty.}\label{fig:muas}
\end{figure}

\begin{algorithm}[!b]
\Fn{\AlCapSty{Multi-Fidelity SRBF for numerical quadrature}}{
Define the initial training sets $\{\mathcal{T}_{\alpha}\}_{\alpha=1}^M$, {with} $\mathcal{T}_{\alpha}=\{{\bf y}_j,G_{\alpha}({\bf y}_j)\}_{j=1}^{\mathcal{J}_{\alpha}}$ \; 
\While(\tcp*[f]{iterate the adaptive sampling process with index $t$}){stopping criteria are not met}{
	Set ${\mathcal{K}^{*,t}_{\alpha}}=\mathcal{J}^t_{\alpha}\,\, \forall \alpha$ \;
	\For{$\alpha=1,\dots,M$}{
		\If(\tcp*[f]{SRBF auto-tuning based on LOOCV}){$\mathcal{J}^t_{\alpha}>5^N$}{
                  Find ${\mathcal{K}^{*,t}_{\alpha}}$ as in Eq. \eqref{eq:loocv} with $\mathcal{C}=[{\mathcal{K}_\alpha^{*,t-1}},
                  {\mathcal{J}^t_\alpha}]$ \;
		} %end If
	} %end For
	Construct the SRBF surrogate $\widehat{G}_1(\mathbf{y})$ as in Eq. \eqref{eq:SRBF} \tcp*{low-fidelity approximation}
	Compute the prediction uncertainty $U_{\widehat{G}_1}(\mathbf{y})$ as in Eq. \eqref{eq:USRBF} \;
	\For(\tcp*[f]{evaluate surrogates of the inter-level errors}){$\alpha=2,\dots,M-1$}{
		Compute the inter-level errors \;
		Construct the SRBF surrogates of the inter-level errors $\tilde{\varepsilon}_{\alpha}$ as in Eq. \eqref{eq:SRBF} \;
		Compute the prediction uncertainty $U_{\tilde{\varepsilon}_{\alpha}}$ as in Eq. \eqref{eq:USRBF} \;
	}
	Construct the multi-fidelity approximation $\widehat{G}_M(\mathbf{y})$ as in Eq. \eqref{eq:ML-MF} \tcp*{MF approximation}
	Compute the multi-fidelity prediction uncertainty $U_{\widehat{G}_M}$ as in Eq. \eqref{eq:ML-MF} \;
	\For(\tcp*[f]{perform parallel infill}){$j=1,\dots,p$}{                          
	     Find ${\bf y}^{\star}={\underset{{\bf y}}{\rm argmax}}[U_{\widehat{G}_M}({\bf y})]$ \;
		 Find $k_j=\mathrm{maxloc}\left[\mathbf{U}(\mathbf{y}^\star)\right]$  as in Eq. \eqref{eq:choice} \;
         Update the training sets $\{\mathcal{T}_{\alpha}\}_{\alpha=1}^{k_j }\cup \{{\bf y}_j^{\star},\widehat{G}_{\alpha}({\bf y}_j^{\star})\}_{\alpha=1}^{k_j}$ \tcp*{considering exact prediction}                                           
         Update the training sets size  $\{\mathcal{J}_{\alpha}^{t+j}\}_{\alpha=1}^{k_j}=\{\mathcal{J}_{\alpha}^t\}_{\alpha=1}^{k_j}+1$ \;
	} %end For          
	\For(\tcp*[f]{perform new simulations}){$j=1,\dots,p$}{                       
                 Evaluate $\{{\bf y}_j^{\star},G_{\alpha}({\bf y}_j^{\star})\}_{\alpha=1}^{k_j}$ \; 
                 Update the training sets $\{\mathcal{T}_{\alpha}\}_{\alpha=1}^{k_j} \cup \{{\bf y}^{\star},G_{\alpha}({\bf y}^{\star})\}_{\alpha=1}^{k_j}$\;                                          
                Update the training sets size $\{\mathcal{J}_{\alpha}^{t+j}\}_{\alpha=1}^{k_j}=\{\mathcal{J}_{\alpha}^t\}_{\alpha=1}^{k_j}+1$\;
           }
           {Get statistics of the quantities of interest by numerical quad. of the surrogate model} \tcp*{UQ}    
           $t=t+1$ \tcp*{Move to the next adaptive sampling iteration}                                              
} % end While
}% end \Fn
\caption{Adaptive multi-fidelity SRBF implementation}
\label{algo:mfm_implementation}
\end{algorithm}

Assuming that the uncertainty associated to the prediction of the lowest-fidelity $U_{\widetilde{G}_1}$ and inter-level errors $U_{\tilde{\varepsilon}_i}$ as uncorrelated, the multi-fidelity approximation $\widehat{G}_M(\mathbf{y})$ of $G_M(\mathbf{y})$ and its uncertainty $U_{\widehat{G}_M}$ read 
\begin{equation}
\label{eq:ML-MF}
\widehat{G}_M(\mathbf{y})\approx\widetilde{G}_1(\mathbf{y})+\sum_{i=1}^{M-1}\tilde{\varepsilon}_i(\mathbf{y})
~~~~~
\mathrm{and}
~~~~~
U_{\widehat{G}_M}(\mathbf{y})=
\sqrt{U^2_{\widetilde{G}_1}(\mathbf{y})+\sum_{i=1}^{M-1}U^2_{{\tilde{\varepsilon}}_i}(\mathbf{y})}.
\end{equation}
%
%The contribution of each fidelity level to $U_{\widehat{G}_M}$ is assessed and used to refine adaptively the training sets as the sampling of the uncertain parameters domain progresses as explained in the following.

Upon having evaluated $U_{\widehat{G}_M}$ the multi-fidelity surrogate is then updated adding a new training point following a two-steps procedure: firstly, the coordinates of the new training point ${\bf y}^\star$ are identified based on the SRBF maximum uncertainty \cite{serani2019-IJCFD}, solving the single-objective maximization problem: 
\begin{eqnarray}\label{eq:MUAS}
{\bf y}^{\star}={\underset{{\bf y}}{\rm argmax}}[U_{\widehat{G}_M}({\bf y})],
\end{eqnarray}
an example (with one fidelity only) is shown in Fig. \ref{fig:muas}. Secondly, once ${\bf y}^\star$ is identified, the training set/sets $\mathcal{T}_{\alpha}$ are updated with a new training point $ \{\mathbf{y}^\star, G_{\alpha}(\mathbf{y}^\star)\}$ with $\alpha=1,\dots,k$, where $k$ is defined as 
\begin{equation}\label{eq:choice}
 k=\mathrm{maxloc}\left[\mathbf{U}(\mathbf{y}^\star)\right]
 ~~~~~
\mathrm{and}
~~~~~
\mathbf{U}(\mathbf{y}^\star)\equiv\{ U_{\widetilde{G}_1}(\mathbf{y}^\star)/\gamma_1, U_{{\hat{\varepsilon}}_1}(\mathbf{y}^\star)/\gamma_2, ...,U_{{\hat{\varepsilon}}_{M-1}}(\mathbf{y}^\star)/\gamma_M \}.
\end{equation}
where $\gamma_{\alpha}$ is the computational cost associated to the $\alpha$-th level.

In the present work, the adaptive sampling procedure starts with five training points (for each fidelity level) located at the domain center and at the lower and upper bounds of each uncertain parameter. SRBF with exact interpolation is enforced for the $\alpha$-th fidelity level until $\mathcal{J}_{\alpha}<5^N$, then the least-squares approximation is used. This because with such a low number of (distant) training points exact interpolation generally provides a better approximation of the quantity of interest than regression. Moreover, the presence of numerical noise of the quantity of interest can not be assessed with few training points, making the regression not useful at the initial stage of the adaptive sampling procedure. Once the least-squares approximation is used, the noise associated to $\alpha$-th fidelity/training set is assessed by Eq. \eqref{eq:rmse} and then normalizing by the range of $G_M$ (to this end, we consider the estimate of the range obtained by the evaluations of $G_M$ available after the first iteration).
Furthermore, to avoid abrupt changes in the SRBF prediction (when regression is used) from one iteration to the next one, the search for the optimal {$\mathcal{K}^*_\alpha$} can be constrained. Herein, defining $t$ as the adaptive sampling iteration,
the problem in Eq. \eqref{eq:loocv} is solved assuming
{$\mathcal{C}=[\mathcal{K}^{*,t-1}_\alpha,\mathcal{J}^t_\alpha]$}, with
{$\mathcal{J}^t_\alpha$} the training set size at the $t$-th iteration, except for the first iteration with least-squares approximation, where no constraint is imposed. 

A deterministic version of the particle swarm optimization algorithm \citep{serani2016-ASC} is used for the solution of the optimization problem in Eq. \eqref{eq:MUAS}.

Since the numerical simulations can be performed with an hardware capable of running $p$ simulations simultaneously, a parallel-infill approach is used in combination with the adaptive sampling method. In the present work, four training points were added to the training sets at each adaptive sampling iteration. The parallel-infill approach computes four sub-iterations for each $t$-th iteration. In each sub-iteration $j$, the training point ${\bf y}^{\star}_{j}$ is identified by Eq. \eqref{eq:choice} and the SRBF prediction $\widehat{G}_{\alpha}({\bf y}^{\star}_{j})$ is added to the $\alpha$-th training set, assuming the SRBF prediction as exact prediction. Once the four training point are identified, the actual simulations are performed and the corresponding outputs are used to update the training sets. Details of the SRBF adaptive sapling procedure are described in Algorithm \ref{algo:mfm_implementation}.

Finally, numerical quadrature is used on the SRBF surrogate model to estimate the expected value of the quantity of interest $\mathbb{E}[\widehat{G}_M]$. {More specifically, $\mathbb{E}[\widehat{G}_M]$
  is approximated using a multi-variate midpoint rule},
with a full-factorial sampling over the SRBF prediction with $S=100^N$ samples, as %
\begin{equation}
\mathbb{E}[G_M] \approx \mathbb{E}[\widehat{G}_M] \approx \frac{1}{S} \sum_{j=1}^S \widehat{G}_M(\bf{y}_j);
\end{equation}
{i.e., a rewriting of Eq. \eqref{eq:MC} where the evaluations of the (discretized) full model $G_{\bm{\alpha}}$ are
  replaced with their multi-fidelity surrogate model counterparts. Sparse-grid quadrature rules can also be considered
  to this end, should a full-factorial sampling be too demanding, see e.g. \cite{gerstner.griebel:adaptive}.}

\section{Problem Formulation and CFD Method}\label{sect:problem_description}
The problem addressed in this manuscript is the forward
UQ analysis of the model-scale resistance ($R_T$) of a RoPax ferry in
straight ahead advancement, subject to two operational uncertainties $\mathbf{y}=[U,T]$, namely the advancement speed ($U$) and the draught ($T$), uniformly distributed within the ranges in Tab. \ref{tab:opcon}.
%
%\begin{equation}
%    \mathbb{E}\left[R_T\right] := \iint R_T(U,T)\rho(U,T)\,\mathrm{d}U\,\mathrm{d}T.
%\end{equation}
%
%Specifically, the uncertain parameters considered $\mathbf{y}=[U,T]$ are the advancement speed ($U$) and the draught ($T$), with uniform PDF.

\begin{figure}[!b]
\centering
\includegraphics[width=0.80\textwidth]{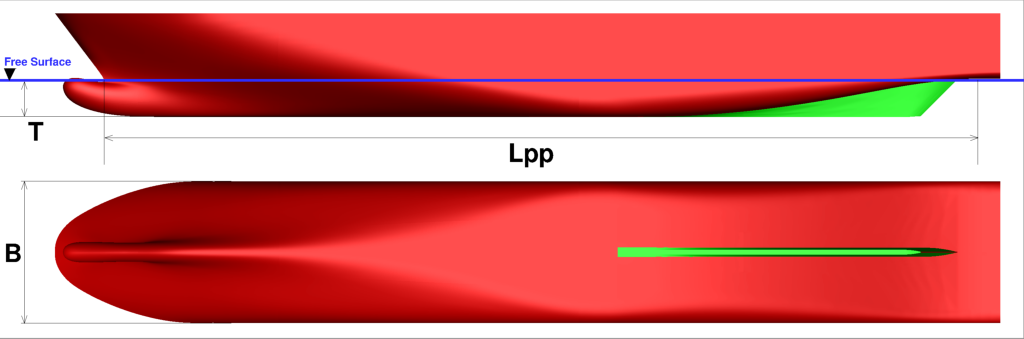}
\caption{RoPax ferry: hull form.}
\label{fig:shape}
\end{figure}

The RoPax ferry is characterized by a length between perpendicular at nominal draught ($L_{\rm PP}$) of
$162.85$~m and a block coefficient $C_B=0.5677$ (see Fig. \ref{fig:shape}).
The parametric geometry of the RoPax is produced with the computer-aided design environment integrated in the CAESES$^{\text{\textregistered}}$ software, developed by FRIENDSHIP SYSTEMS AG, and made available in the framework of the H2020 EU Project Holiship.
%; fore perpendicular is defined as the intersection with the unperturbed free surface level, whereas, the aft perpendicular corresponds to the rudder axys.
The analysis is performed at model scale with a scale factor equal to $27.14$. The main dimensions and the operative conditions are summarized in Tab. \ref{tab:opcon}. The advancement speed ranges from $12$ to $26$ knots at full scale and the draught variation is $\pm 10\%$ of the nominal draught, which corresponds to a variation of about $\pm 15\%$ of the nominal displacement. The corresponding range in Froude number ${\rm Fr}={U}/{\sqrt{g L_{\rm PP}}}$ is $[0.154,0.335]$, whereas the variation in Reynolds number (at model scale) ${\rm Re}=\rho U L_{\rm PP}/\mu=U L_{\rm PP}/\nu$ is $[6.423\cdot10^6,1.392\cdot10^7]$, where $\rho=998.2$ kg/m$^3$ is the water density, $\nu=\mu/\rho=1.105\cdot10^{-6}$ m$^2$/s the kinematic viscosity and $g=9.81$ m/s the gravitational acceleration.

\begin{table}[!b]
	\begin{center}
		\caption{Main geometrical details and operative conditions of the RoPax ferry (model scale $1:27.14$).}
		\label{tab:opcon}
		\begin{tabular}{lcccc}
			\toprule
			\textbf{Description} & \textbf{Symbol} & \textbf{Full Scale} & \textbf{Model Scale} & \textbf{Unit}   \\
			\midrule
			Length between perpendiculars & $L_{\rm PP}$ & $162.85$          & $6.0$              & m     \\
			Beam                          & $B$          & $29.84$           & $1.0993$             & m    \\
			Block coefficient             & $C_B$        & $0.5677$          & $0.5677$           & --   \\
			Nominal displacement          & $\nabla$     & $19584.04$        & $0.9996$           & m$^3$ \\
			Nominal draught               & $T_n$        & $7.10$            & $0.261660$ & m \\
			Draught range                 & $T$  & $[7.812, 6.391]$ & $[0.2355, 0.2878]$ & m \\
			Speed range                   & $U$  & $[6.173, 13.376]$ & $[1.185, 2.567]$ & m/s   \\
			Froude range                  & Fr & $[0.154, 0.335]$ & $[0.154, 0.335]$ & --   \\
			Reynolds range                & Re & $[9.081\cdot10^8, 1.968\cdot10^9]$ & $[6.423\cdot 10^6, 1.392\cdot 10^7]$ & --   \\
			\bottomrule
		\end{tabular}                                           
	\end{center}
\end{table}
%Draught range                 & $T_{\rm MIN}\div T_{\rm MAX}$ & $7.812\div6.391$ & $0.235494\div0.287826$ & m \\
%Speed range                   & $U_{\rm MIN}\div U_{\rm MAX}$ & $6.173\div13.376$ & $1.185\div2.567$ & m/s   \\

The hydrodynamics performance of the RoPax ferry {for fixed values of the uncertain parameters} is assessed by the RANS code $\chi$navis developed at CNR-INM \cite{dimascio2007,dimascio2009,broglia2018}. It is based on a finite volume {discretization of the RANS equations},
with variables collocated at the cell centers. Turbulent stresses are taken into account by the Boussinesq hypothesis, with the Spalart-Allmaras turbulence model. Free-surface effects are taken into account by a single-phase level-set algorithm. Wall-functions are not adopted, therefore the wall distance $y^+=1$ is ensured on the wall.

\begin{figure}[!t]
\centering
\subfigure[Boundary conditions]{\includegraphics[width=0.45\textwidth]{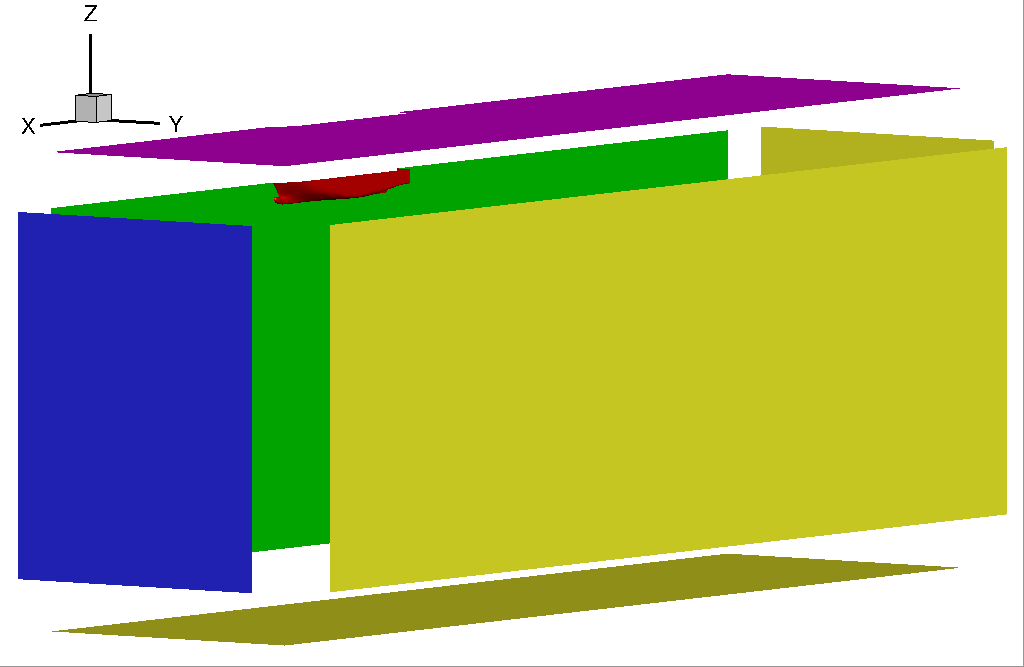}}
\subfigure[Computational mesh] {\includegraphics[width=0.45\textwidth]{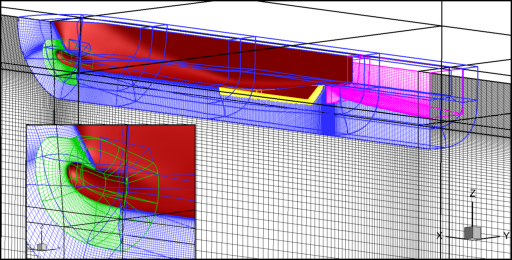}}
\caption{Numerical setup: boundary conditions and computational grid.}
\label{fig:bc}
\end{figure}

The computational domain extends to 2$L_{\rm PP}$ in front of the hull, 3$L_{\rm PP}$ behind, and 1.5$L_{\rm PP}$ {on} each side; a depth of 2$L_{\rm PP}$ is imposed (see Fig. \ref{fig:bc}). On the solid walls (in red in the figure), the velocity is set equal to zero, whereas zero normal gradient is enforced on the pressure field; at the (fictitious) inflow boundary (blue), the velocity is set to the undisturbed flow value and the pressure is extrapolated from inside; the dynamic pressure is set to zero at the outflow (yellow), whereas the velocity is extrapolated from inner points. On the top boundary, which remains always in the air region, fluid dynamic quantities are extrapolated from inside (purple).
Taking advantage of the symmetry of the flow relative to the $y=0$ plane, computations are performed for half ship only, {and the} usual symmetry boundary conditions are enforced on the symmetry longitudinal plane (in green in Fig.~\ref{fig:bc}{\it a}).

\begin{figure}[!b]
\centering
\includegraphics[width=1\textwidth]{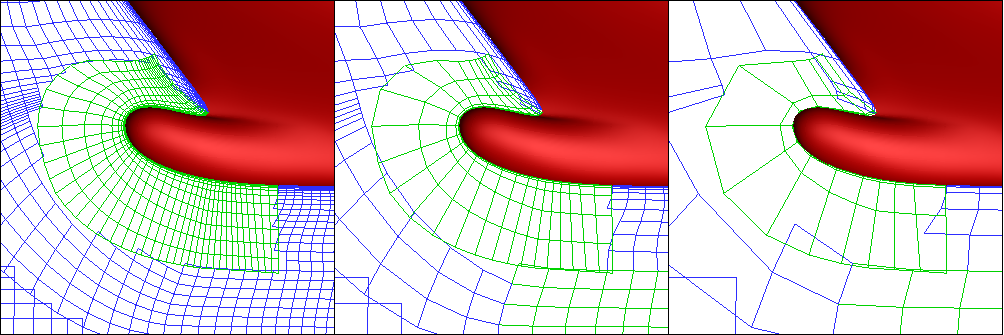}
\caption{RoPax grids: detail of the bow region, left to right: $\mathcal{M}_3$ (699k), $\mathcal{M}_2$ (87k), $\mathcal{M}_1$ (11k).}
\label{fig:grids}
\end{figure}
\begin{figure}[!b]
\centering
\includegraphics[width=1\textwidth]{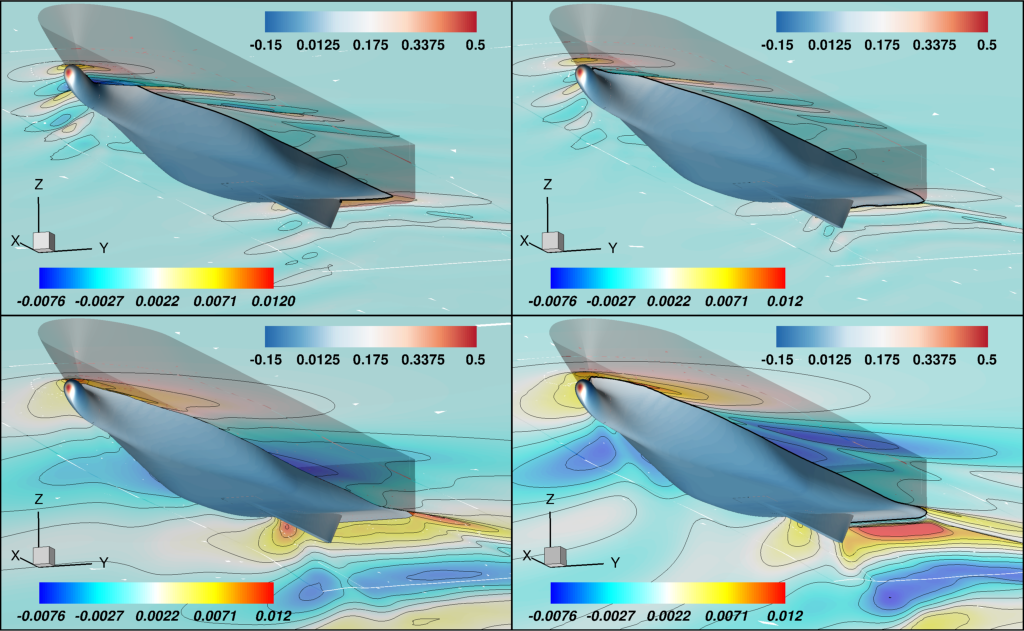}
\caption{RoPax ferry, $\mathcal{M}_4$ CFD results in terms of non-dimensional wave pattern and surface pressure for:
	${\rm Fr}=0.193$, $T=3.9249\cdot10^{-2}L_{\rm PP}$ and $T=4.7971\cdot10^{-2}L_{\rm PP}$, top row left and right;
	${\rm Fr}=0.335$, $T=3.9249\cdot10^{-2}L_{\rm PP}$ and $T=4.7971\cdot10^{-2}L_{\rm PP}$, bottom row left and right.}
	\label{fig:RoPax1}
\end{figure}

\begin{figure}[!t]
\centering
\includegraphics[width=1\textwidth]{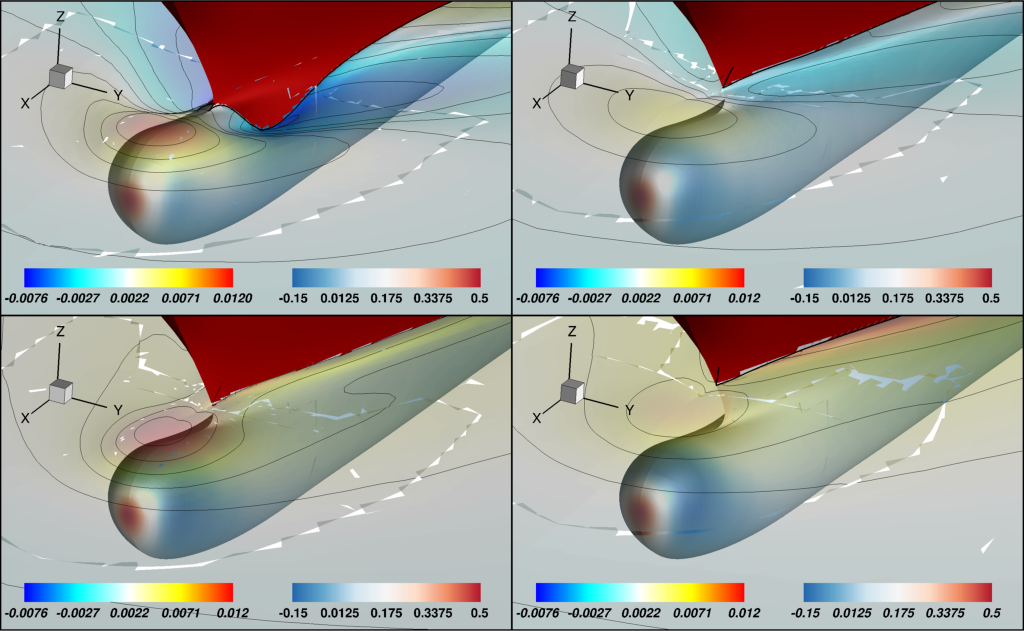}
\caption{RoPax ferry, enlarged view of the bow region as in Fig. \ref{fig:RoPax1}.}
\label{fig:RoPax2}
\end{figure}

The computational grid is composed by 60 adjacent and partially overlapped blocks; Fig. \ref{fig:bc}{\it b} shows a particular of the block structures in the region around the ship hull and the computational mesh on the symmetry plane. Taking the advantage of a Chimera overlapping approach, the meshes around the skeg and around the bow are generated separately from the mesh around the hull; a background Cartesian grid is then built and the whole grid is assembled by means of an in-house overlapping grid pre-processor. The final mesh counts for a total of about 5.5M of control volumes for half the domain. The numerical solutions are computed by means of a full multi-grid--full approximation scheme (FMG--FAS), with four grid levels (from coarser to finer: $\mathcal{M}_1$, $\mathcal{M}_2$, $\mathcal{M}_3$, and $\mathcal{M}_4$), each obtained from the next finer grid with a coarsening ratio equal to 2, along each curvilinear direction. In the FMG--FAS approximation procedure, the solution is computed on the coarsest grid level first. Secondly, it is approximated on the next finer grid and the solution is iterated by exploiting all the coarser grid levels available with a V-Cycle. The process is repeated up to the finest grid level. For the present UQ problem all the four grid levels are used; to note, the number of grid volumes ranges from 5.5M for the finest mesh, down to 11K for the coarsest one. To provide an idea about the different mesh resolution between the grid levels, Fig. \ref{fig:grids} shows a particular of the grid at the bow region for $\mathcal{M}_3$, $\mathcal{M}_2$ and $\mathcal{M}_1$ grids ($\mathcal{M}_4$ is in the insert of Fig. \ref{fig:bc}{\it b}).

Based on the grid refinement ratio chosen, a normalized computational cost for the $\alpha$-th grid level can be estimated as:
\begin{equation}
  \label{eq:cost}
  \text{cost}({\alpha}) = 8^{\alpha-1}
\end{equation}
with $\alpha=1,\dots,4$.
In the FMG-FAS scheme the computation on the $\alpha$-th grid level involves computations on all the coarser meshes $\mathcal{M}_1\dots\mathcal{M}_{\alpha-1}$. However, with the estimation in Eq. \eqref{eq:cost}, only the cost of the highest-fidelity level samples is taken into account, i.e. the computations on the coarser grids are considered negligible.

Fig. \ref{fig:RoPax1} shows an overview of the numerical solutions obtained for different conditions in terms of wave pattern and hull surface; wave height (as elevation with respect to the unperturbed level) and surface pressure are reported in non-dimensional values, making the height non dimensional with $L_{\rm PP}$ and the pressure with $\rho U^2 L_{\rm PP}^2$, as usual. %Figures refer to the finest mesh level (i.e. $G_4$).
A clear, and obvious, Froude number dependency is seen for the wave patterns; at the lower speed shown, the free surface is weakly perturbed (to note, the same color range has been used for all the panels), whereas, at higher Froude, a clear kelvin pattern is seen. Also, at higher speed, the formation of a well defined transom wave system is observed, including the presence of the classical rooster tail.
It is also worth to observe the influence of the draught on the wave system; in particular at the lower speed reported, and the smaller draught, the rear part of the bulbous is partially dry (better seen in the enlarged views reported in Fig.~\ref{fig:RoPax2}). The region of very low pressure, caused by the flow acceleration around the bow, is obviously the cause. For all cases, the high pressure in the stagnation point at the bow prevents the bow to be outside the water, as it is at the rest conditions at least for the nominal and the smaller draughts (see Fig.~\ref{fig:shape}). For higher speed, the larger draught condition causes a stronger rooster tail system at the stern, with higher crest and trough.

\section{Numerical Results} \label{sect:numerical_results}
\begin{figure}[!t]
	\centering
	\subfigure[Expected value]{\includegraphics[width=1\textwidth]{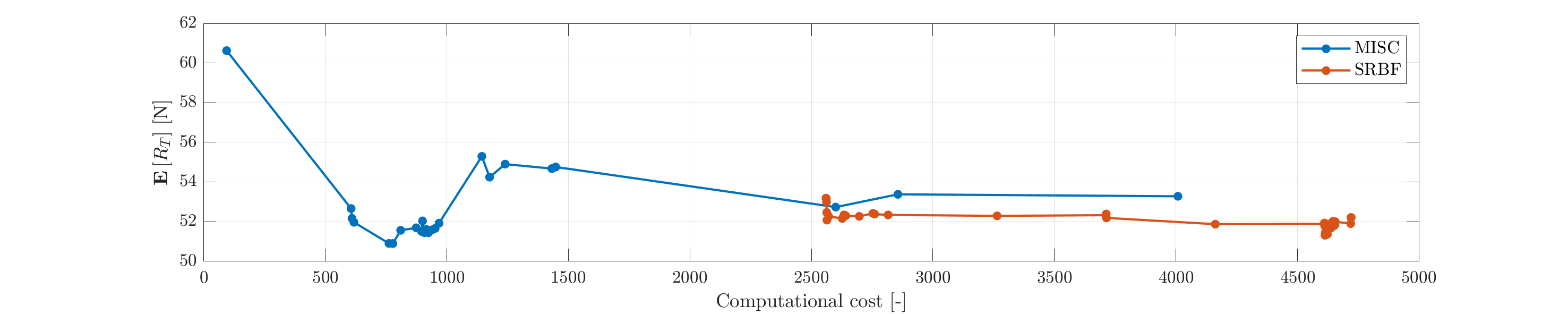}}
	\subfigure[Standard deviation]{\includegraphics[width=1\textwidth]{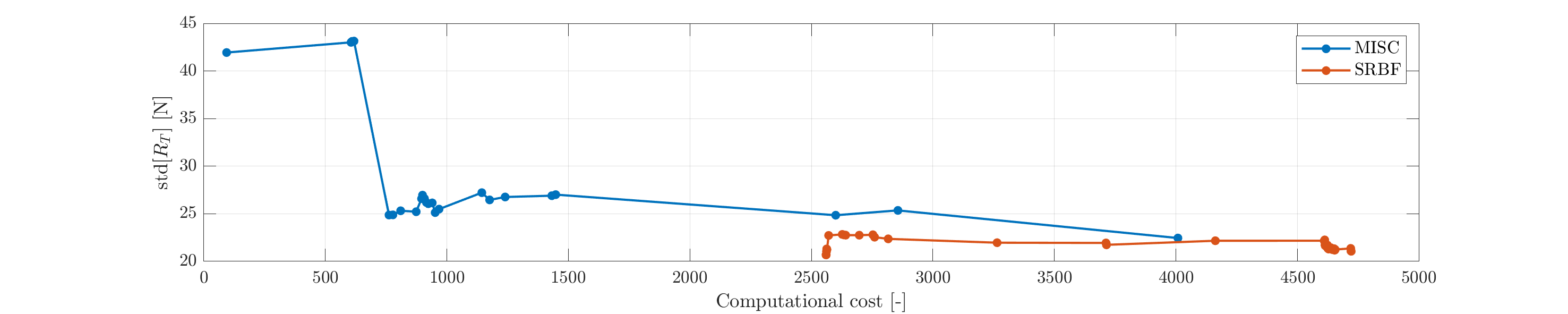}}
	\caption{Comparison of MISC and SRBF results: convergence of the expected value and standard deviation of $R_T$ versus computational cost. For the sake of readability, we plot the results for MISC starting with the 4th iteration. The first three iterations correspond to computational cost 1, 12, and 14, and give poor results.}
	\label{fig:EV_comparison}
\end{figure}

In this section the behavior of the two methods applied to the problem just described is illustrated.
To fairly compare the outcomes of the two UQ methodologies presented, no termination criteria are used for both methods, that are stopped when sufficient and comparable computational costs have been reached.

\begin{table}[!b]
\centering
\caption{Comparison of MISC and SRBF results: computational cost, expected value and standard deviation of $R_T$ for an intermediate iteration and at the final iteration.}\label{table:comparison}
\begin{tabular}{ccccc}
\toprule
\bf UQ method & \multicolumn{2}{c}{MISC} & \multicolumn{2}{c}{SRBF}\\
\midrule
Iteration                    & 14      & final   & 8       & final \\
Computational cost           & 900     & 4008    & 4638    & 4721 \\
$\mathbb{E}\left[R_T\right]$ & 52.0    & 53.3    & 51.7    & 52.2 \\
$\text{std}\left[R_T\right]$ & 26.9    & 22.4    & 21.4    & 21.0\\
\bottomrule
\end{tabular}
\end{table} 

\begin{figure}[!b]
	\centering
	\subfigure[MISC, iteration 14]{\includegraphics[width=0.45\textwidth]{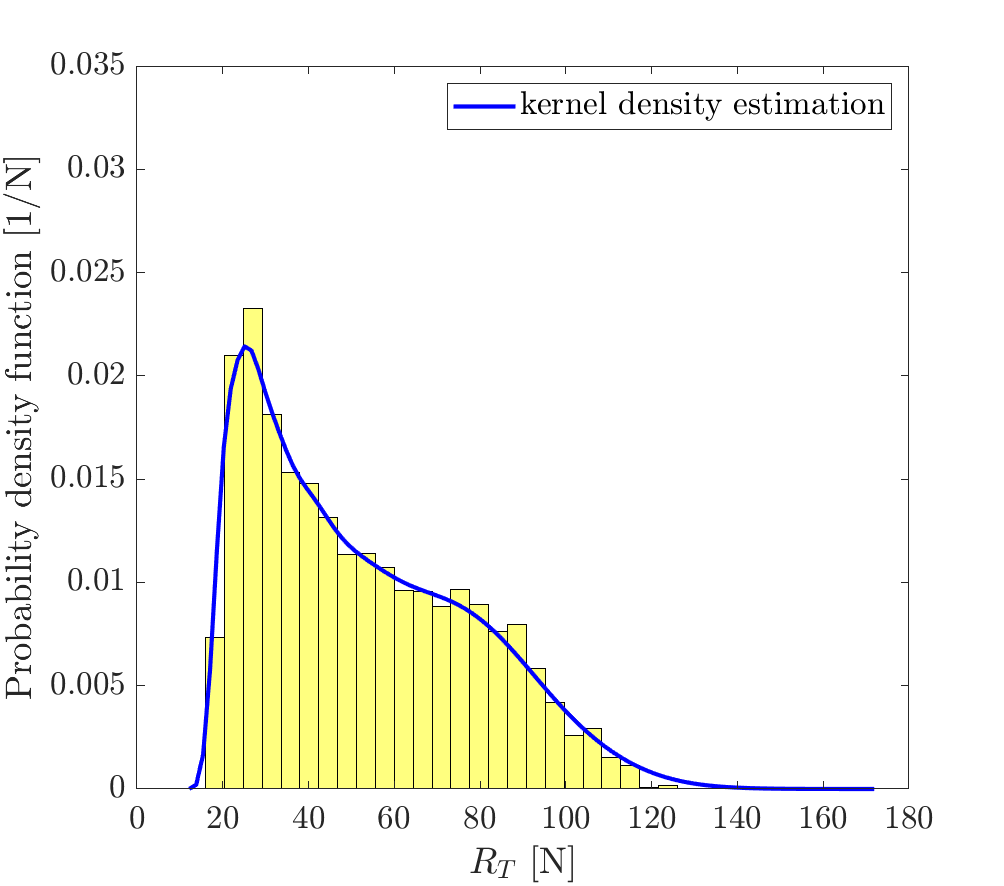}}
	\subfigure[MISC, final iteration]{\includegraphics[width=0.45\textwidth]{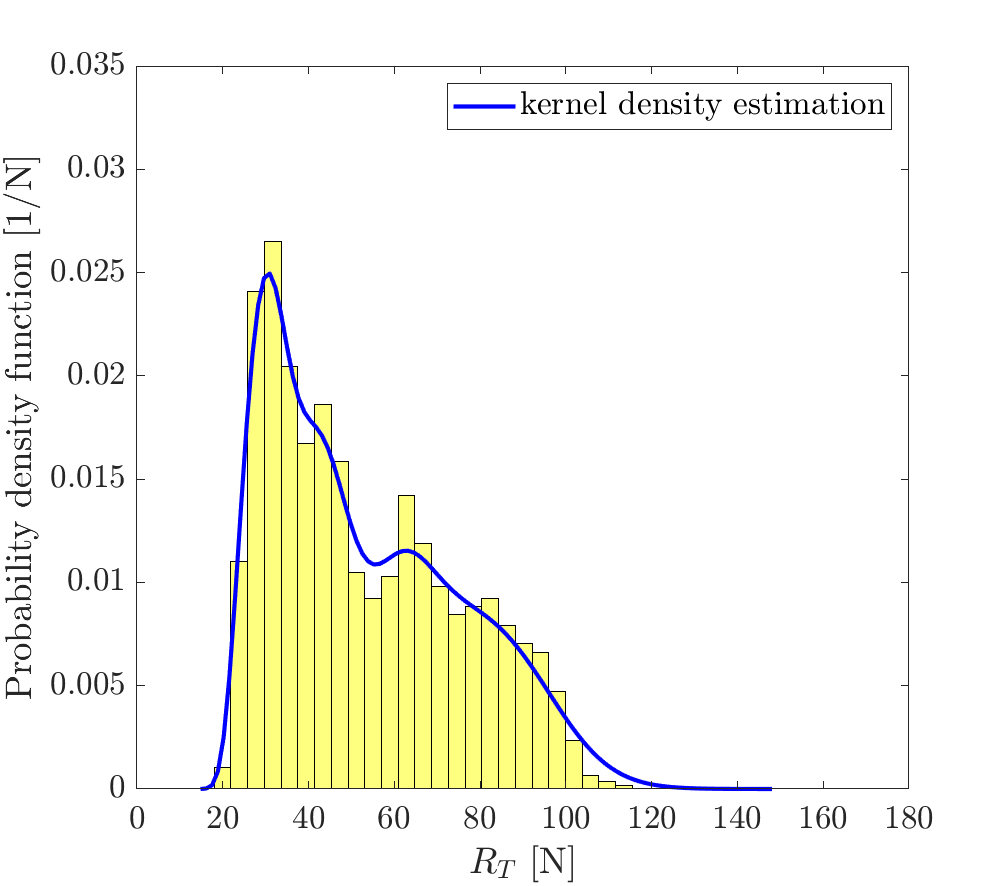}}
	\subfigure[SRBF, iteration 8]{\includegraphics[width=0.45\textwidth]{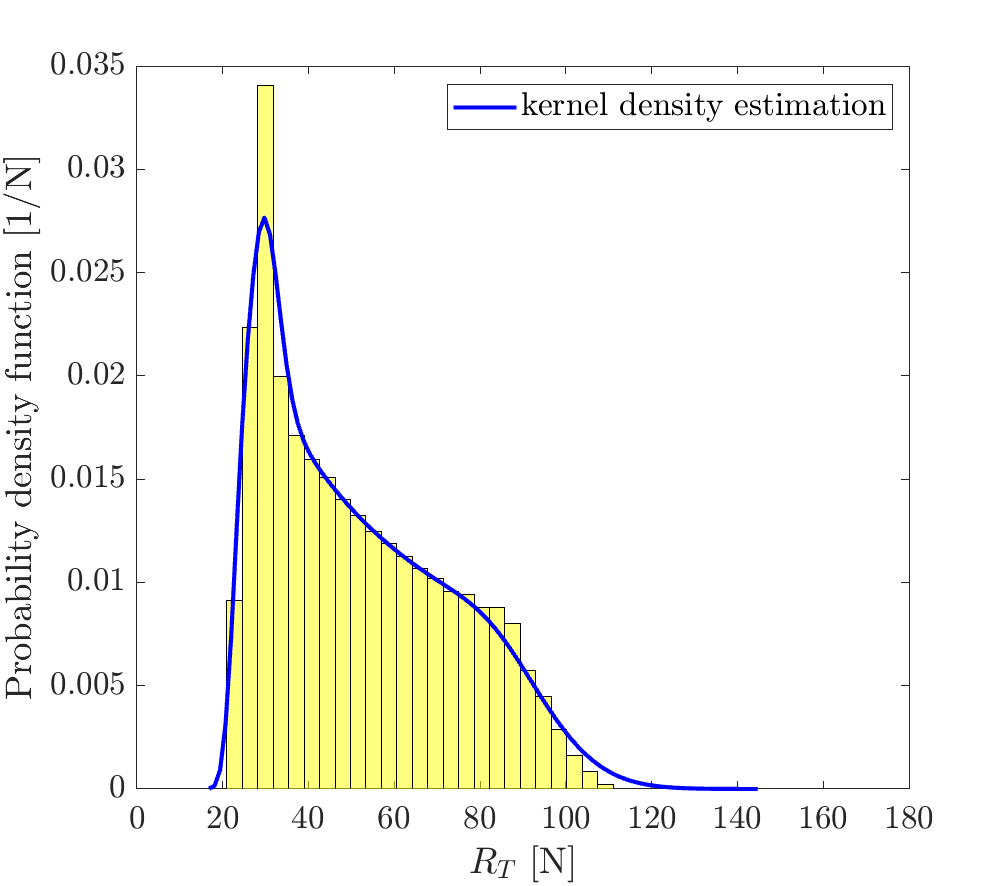}}
	\subfigure[SRBF, final iteration]{\includegraphics[width=0.45\textwidth]{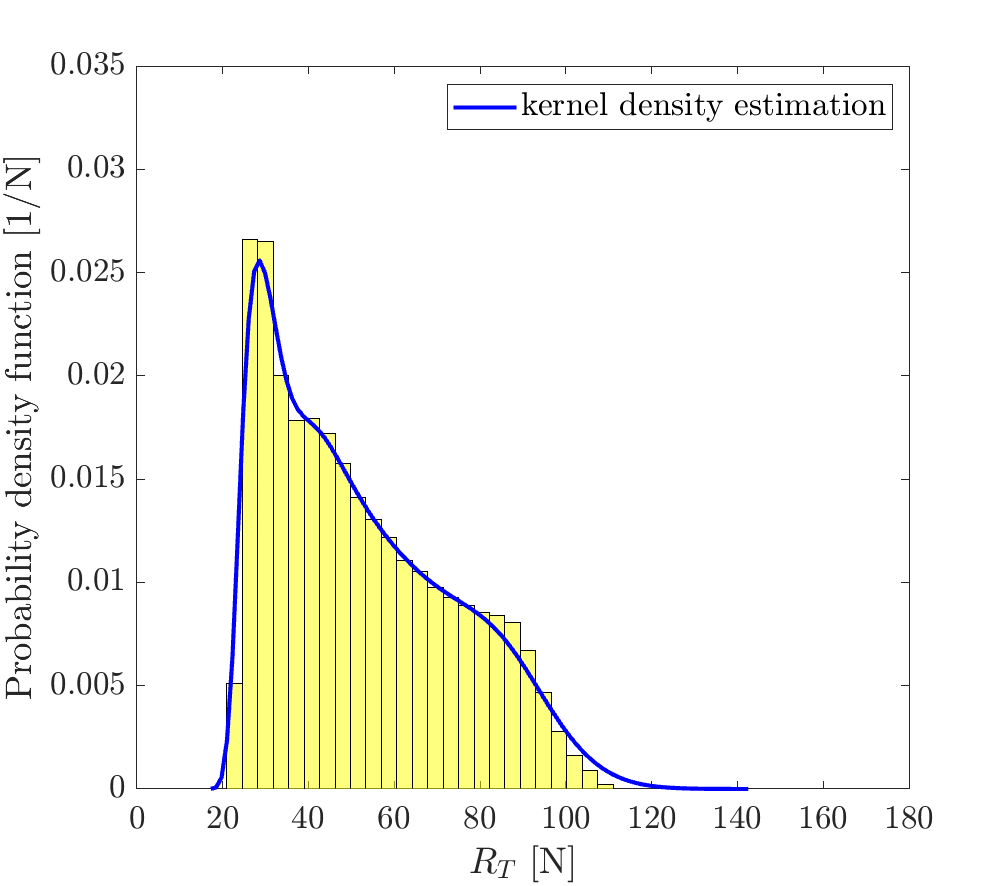}}
	\caption{Comparison of MISC and SRBF results: density function and histogram of (rescaled) frequencies of $R_T$. The histogram is computed using 25 bins for $n=10000$ samples of the response surface. The probability density function is estimated using the kernel smoothing estimate given by the Matlab's ksdensity function, enforcing positive support and with automatic selection of the bandwidth.
%the Matlab kernel smoothing function estimate given by \texttt{ksdensity} with positive support and automatic choice of the bandwidth evaluated on 100 equally-spaced points.
}
	\label{fig:pdf_comparison}
\end{figure}

\begin{figure}[!b]
	\centering
%	\subfigure[MISC, iteration 14]{\includegraphics[width=0.24\textwidth]{Fig_MISC/MISC_Contour_iter13_50by50.png}}
%	\subfigure[MISC, final iteration]{\includegraphics[width=0.24\textwidth]{Fig_MISC/MISC_Contour_iter30_50by50.png}}
%          \subfigure[SRBF, iteration 8]{\includegraphics[width=0.24\textwidth]{Fig_RBF/Contour-Iter0033.png}}
%	\subfigure[SRBF, final iteration]{\includegraphics[width=0.24\textwidth]{Fig_RBF/Contour-Iter0053.png}}
	\subfigure[MISC, iteration 14]{\includegraphics[width=0.45\textwidth]{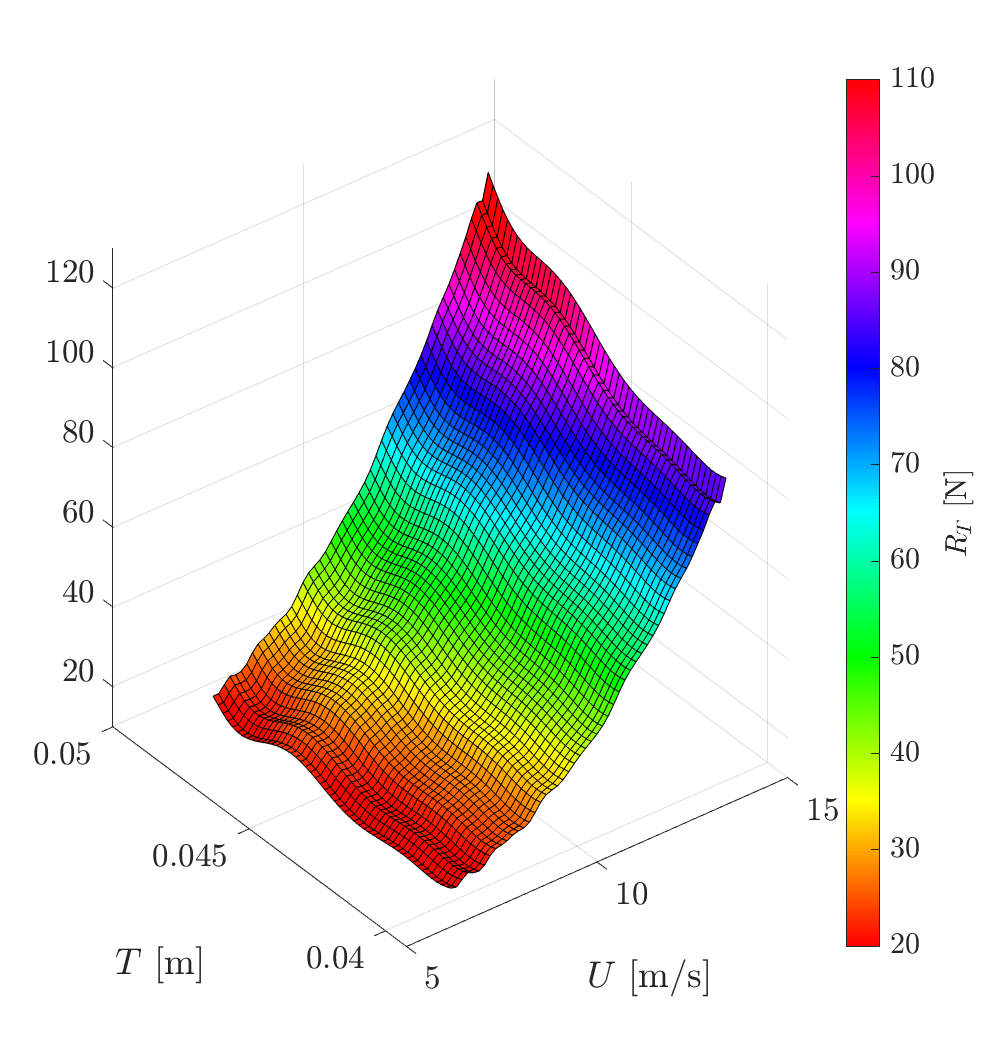}}
	\subfigure[MISC, final iteration]{\includegraphics[width=0.45\textwidth]{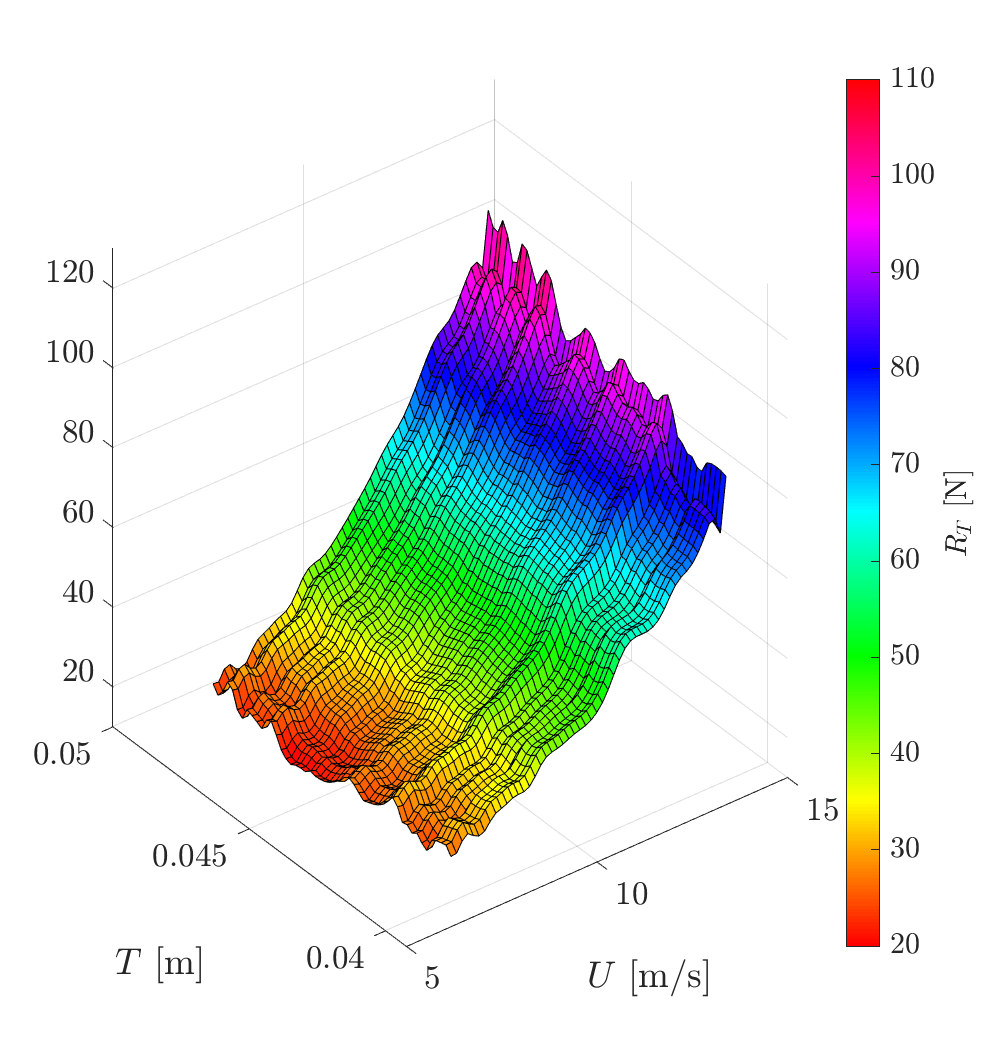}}
	\subfigure[SRBF, iteration 8]{\includegraphics[width=0.45\textwidth]{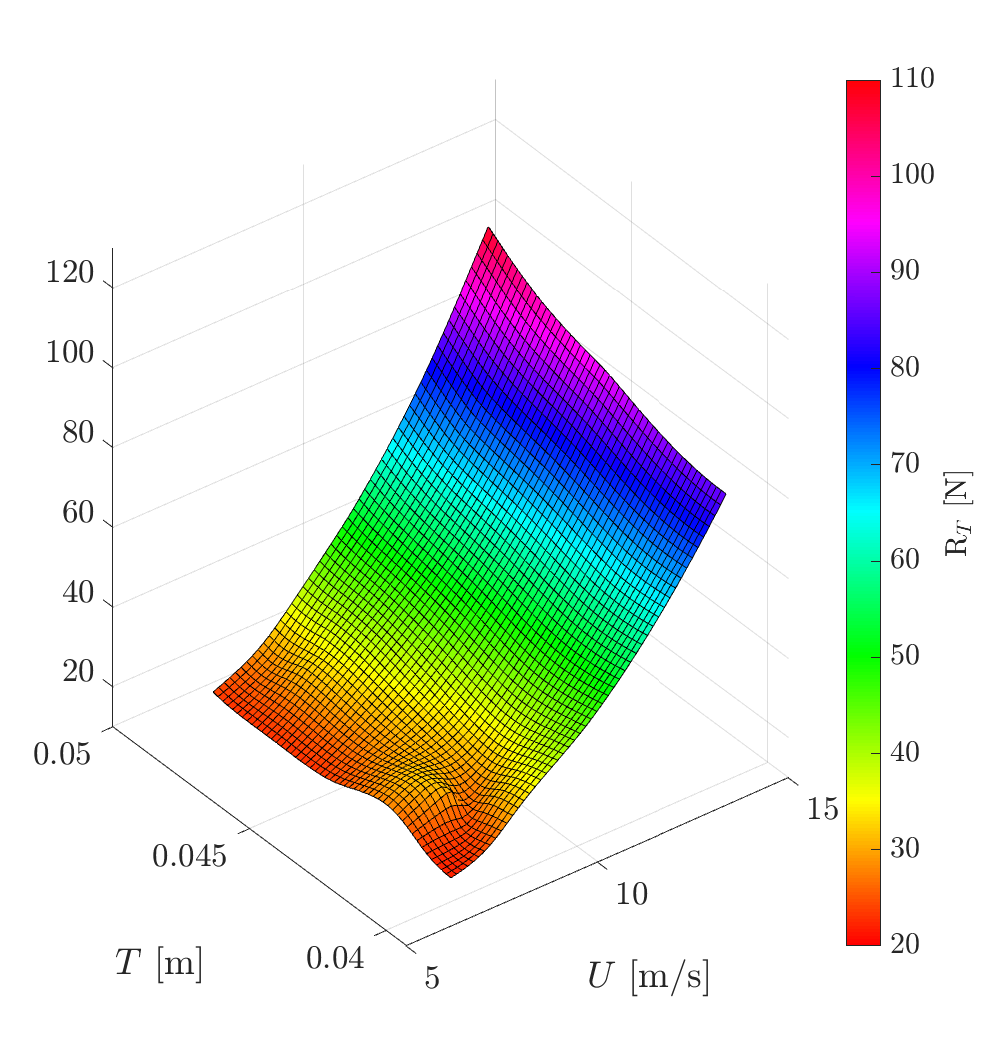}}
	\subfigure[SRBF, final iteration]{\includegraphics[width=0.45\textwidth]{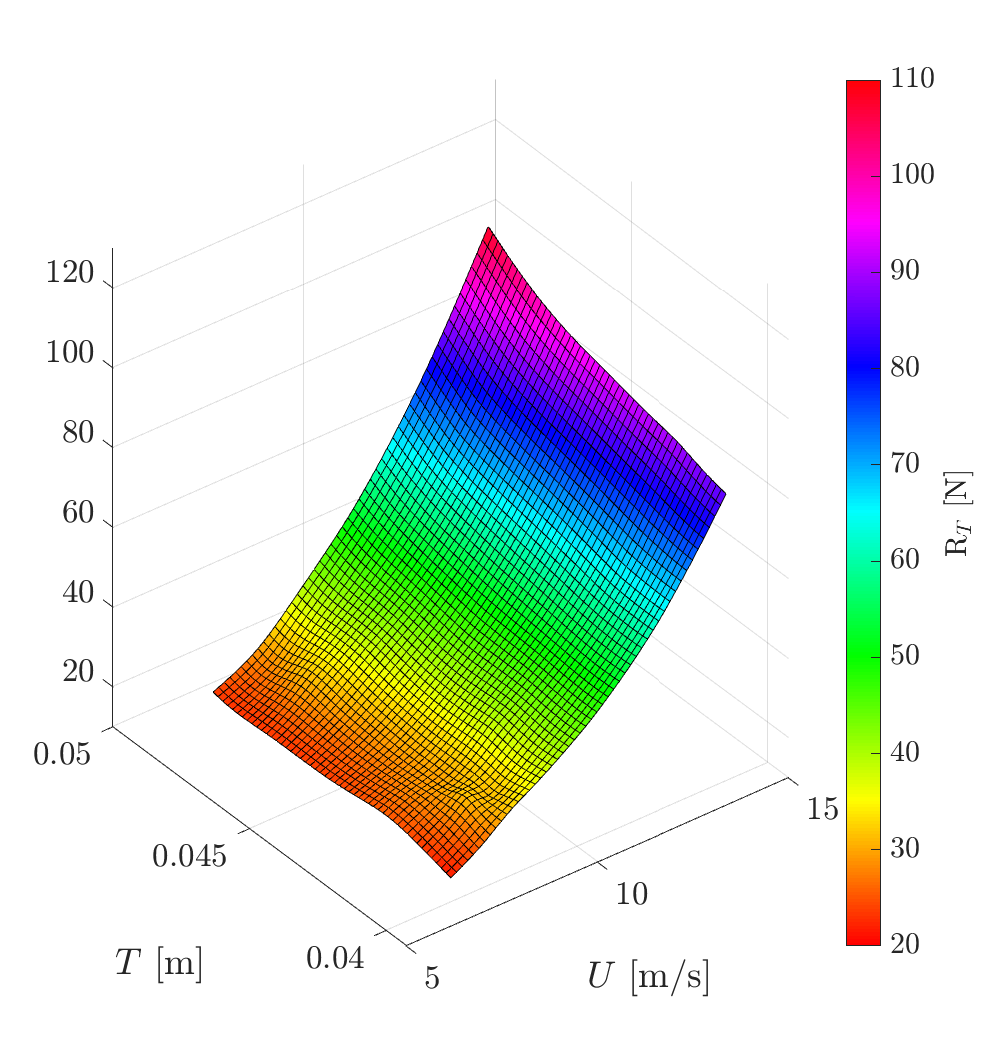}}
	\caption{Comparison of MISC and SRBF results: response surfaces.}
	\label{fig:response_surface_comparison}
\end{figure}

\begin{figure}[!b]
	\centering
	\subfigure[MISC, iteration 14]{\includegraphics[width=0.45\textwidth]{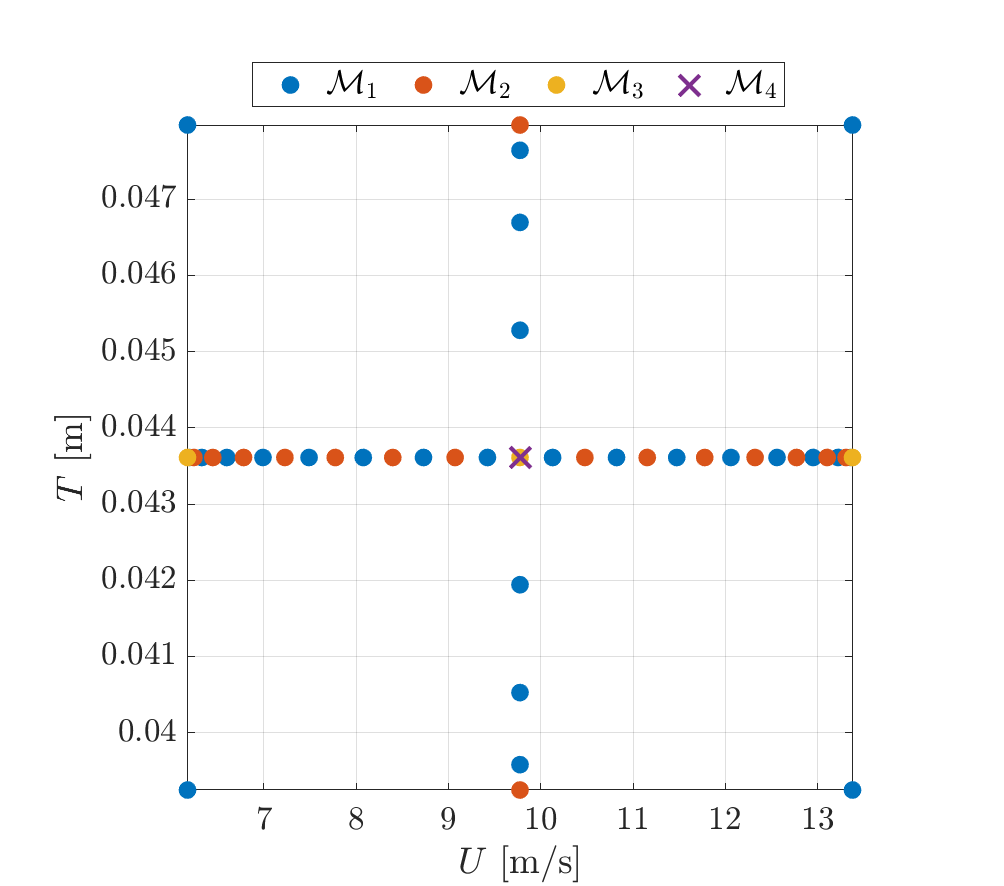}}
	\subfigure[MISC, final iteration]{\includegraphics[width=0.45\textwidth]{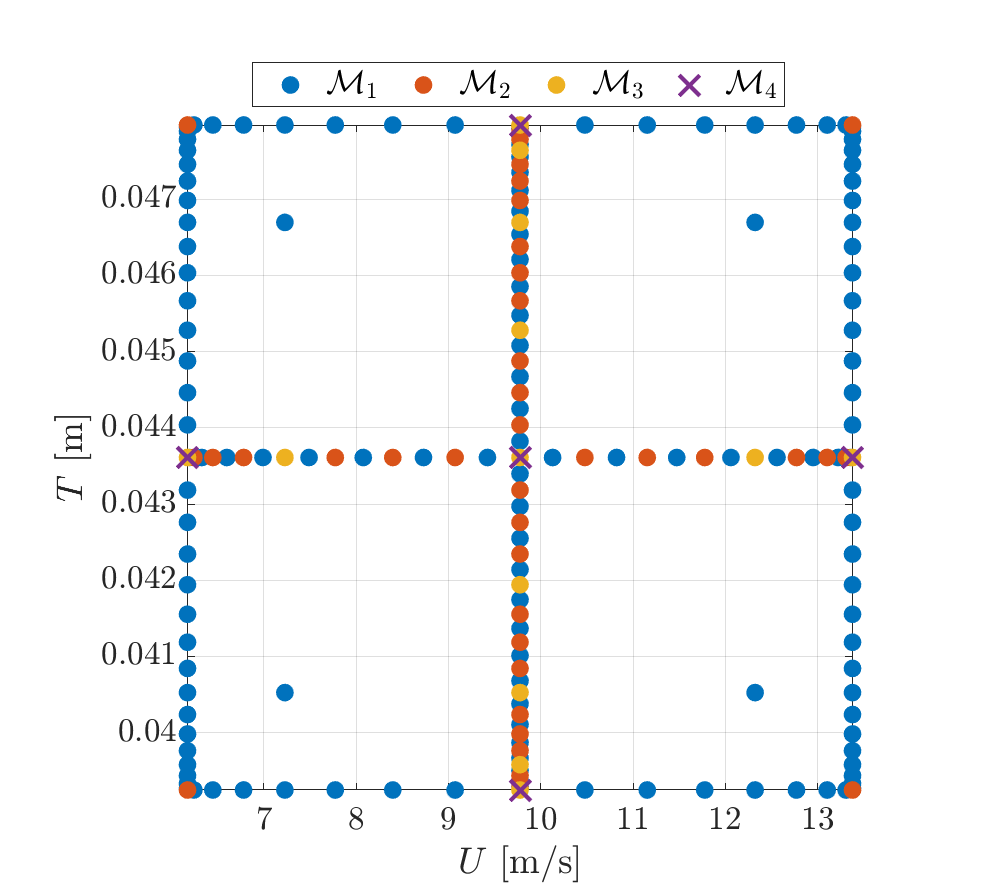}}
           \subfigure[SRBF, iteration 8]{\includegraphics[width=0.45\textwidth]{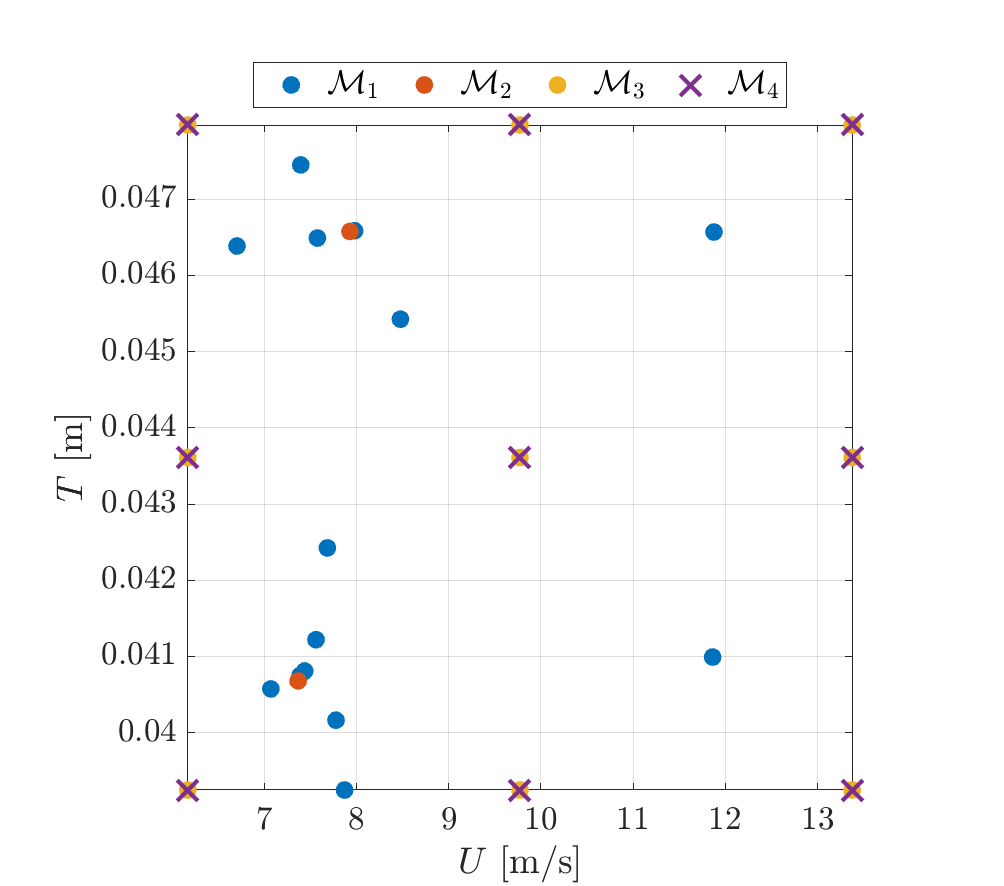}}
	\subfigure[SRBF, final iteration]{\includegraphics[width=0.45\textwidth]{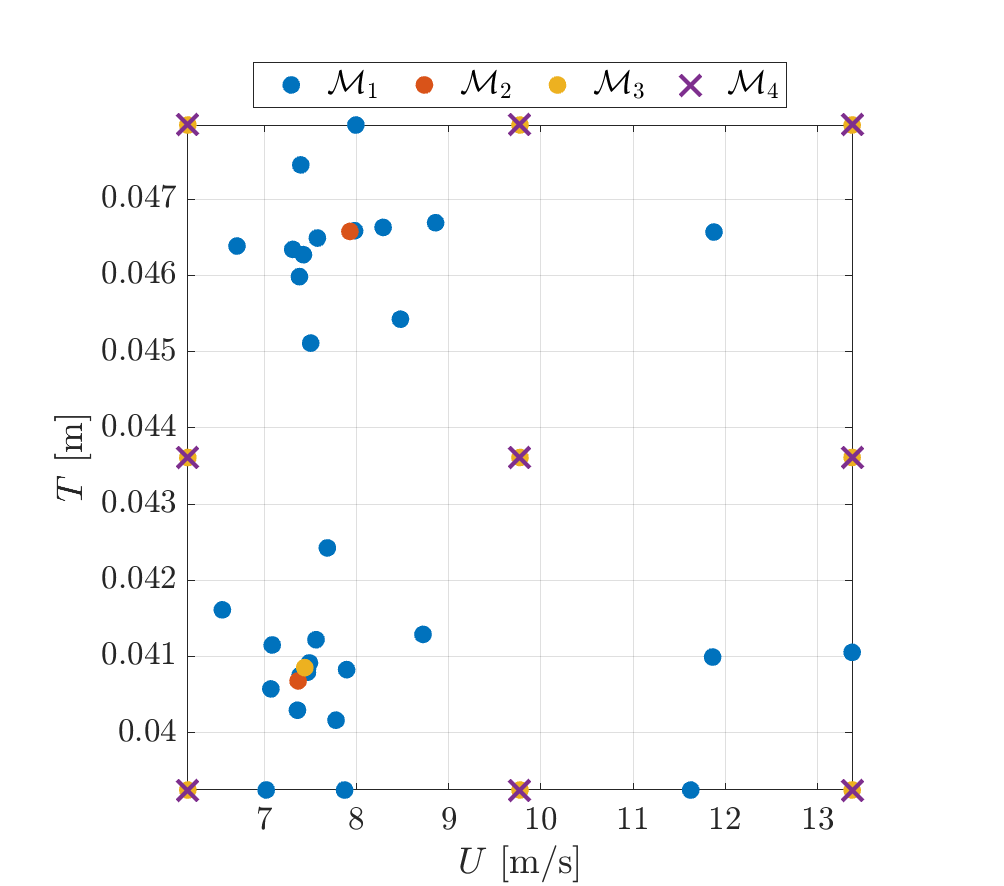}}
	\caption{Comparison of MISC and SRBF results: points in parameter space. Note that every point required on grid $\mathcal{M}_i$ is required also on all the grids with lower refinement level.}
	\label{fig:points_comparison}
\end{figure}

\begin{figure}[!b]
	\centering
	\subfigure[MISC, cumulative number of simulations]{\includegraphics[width=0.45\textwidth]{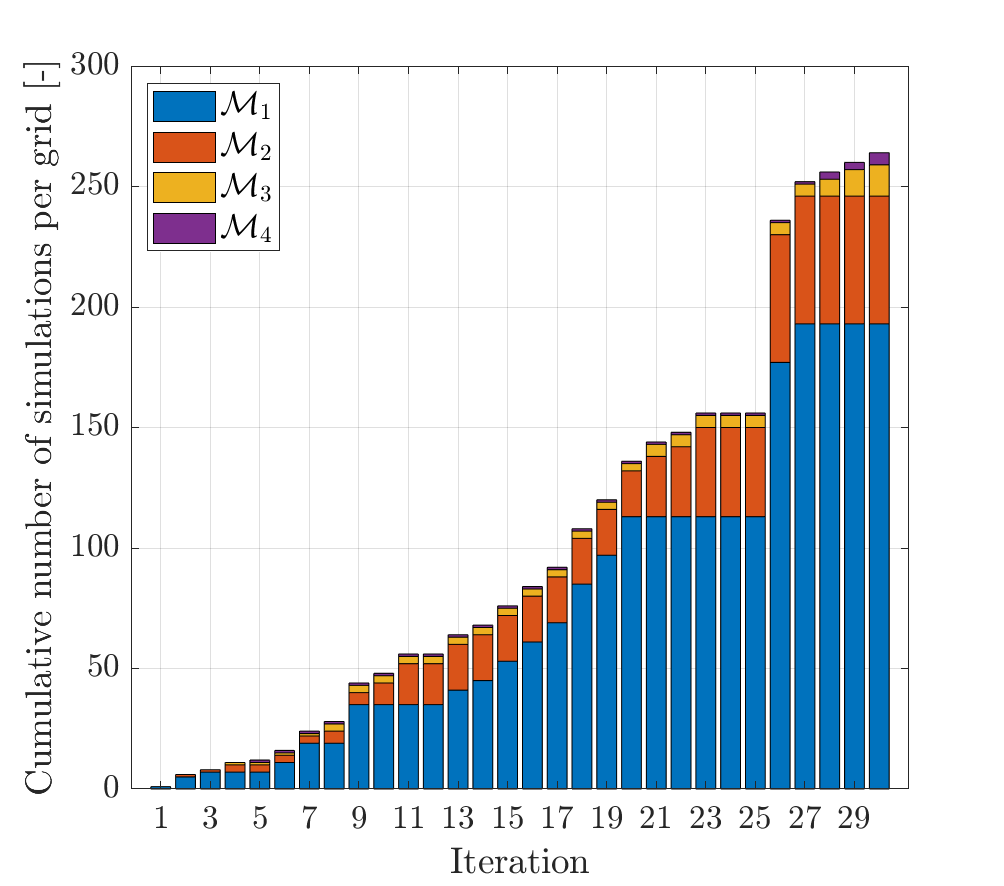}}
	\subfigure[MISC, number of new simulations]{\includegraphics[width=0.45\textwidth]{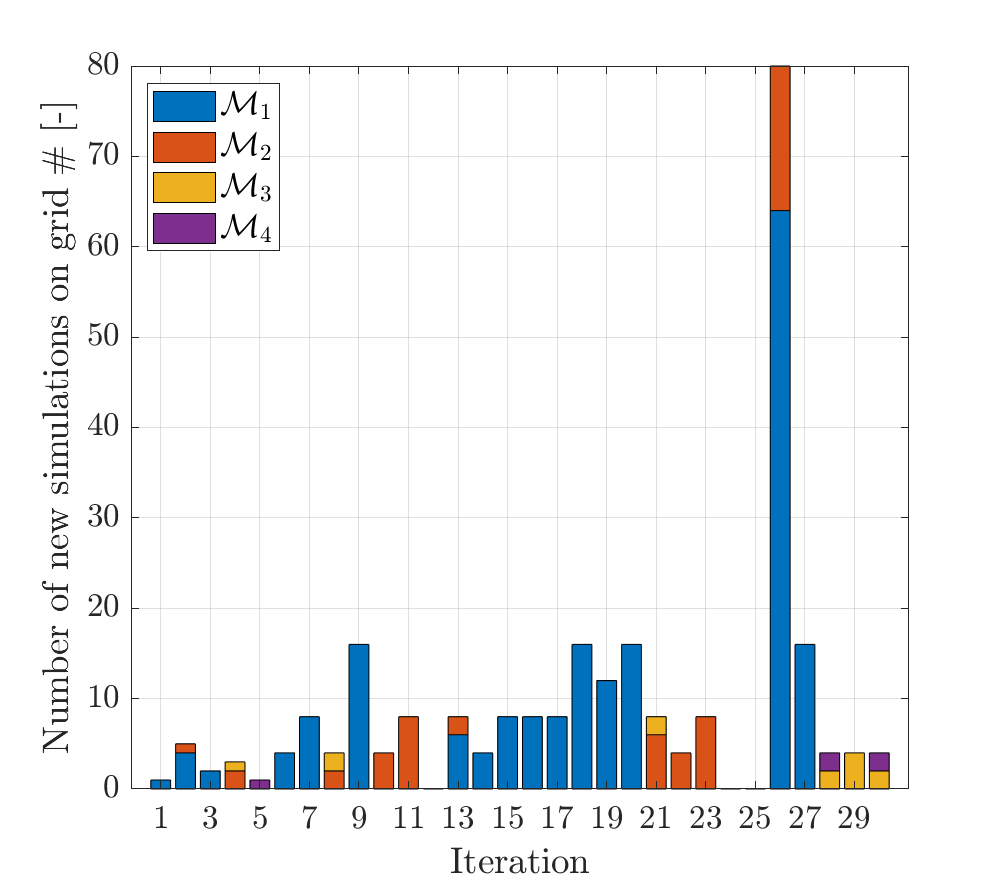}}
           \subfigure[SRBF, cumulative number of simulations]{\includegraphics[width=0.45\textwidth]{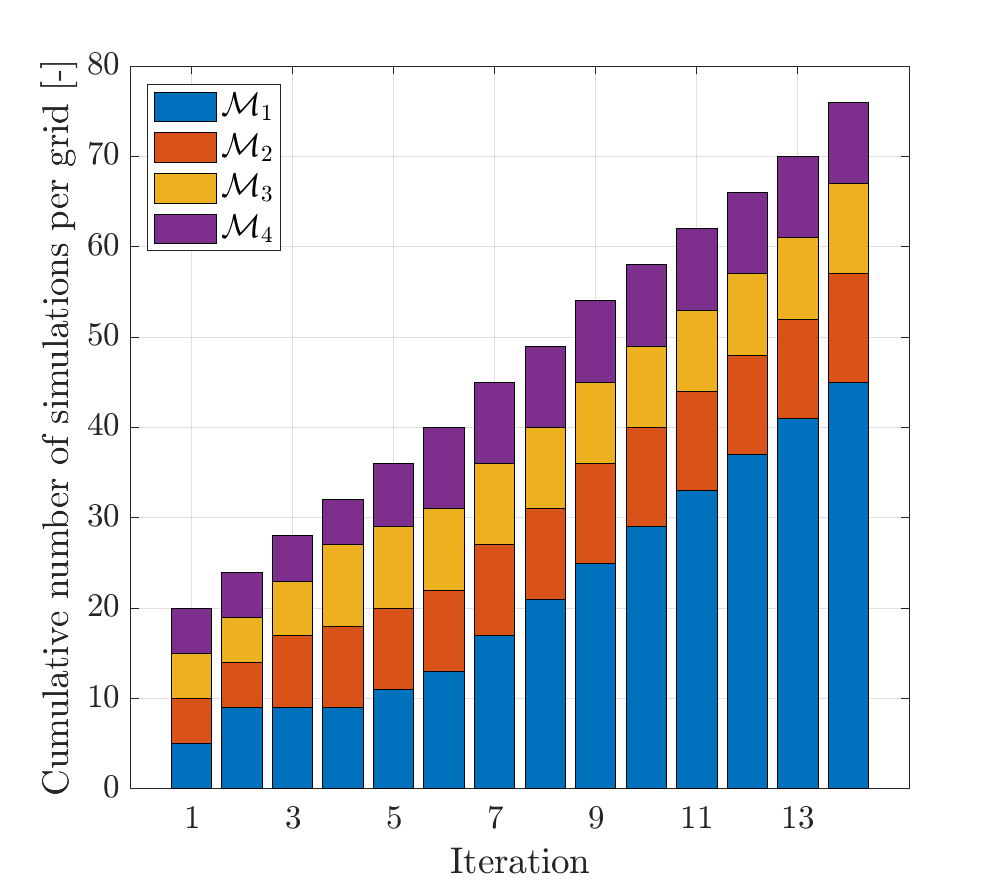}}
           \subfigure[SRBF, number of new simulations]{\includegraphics[width=0.45\textwidth]{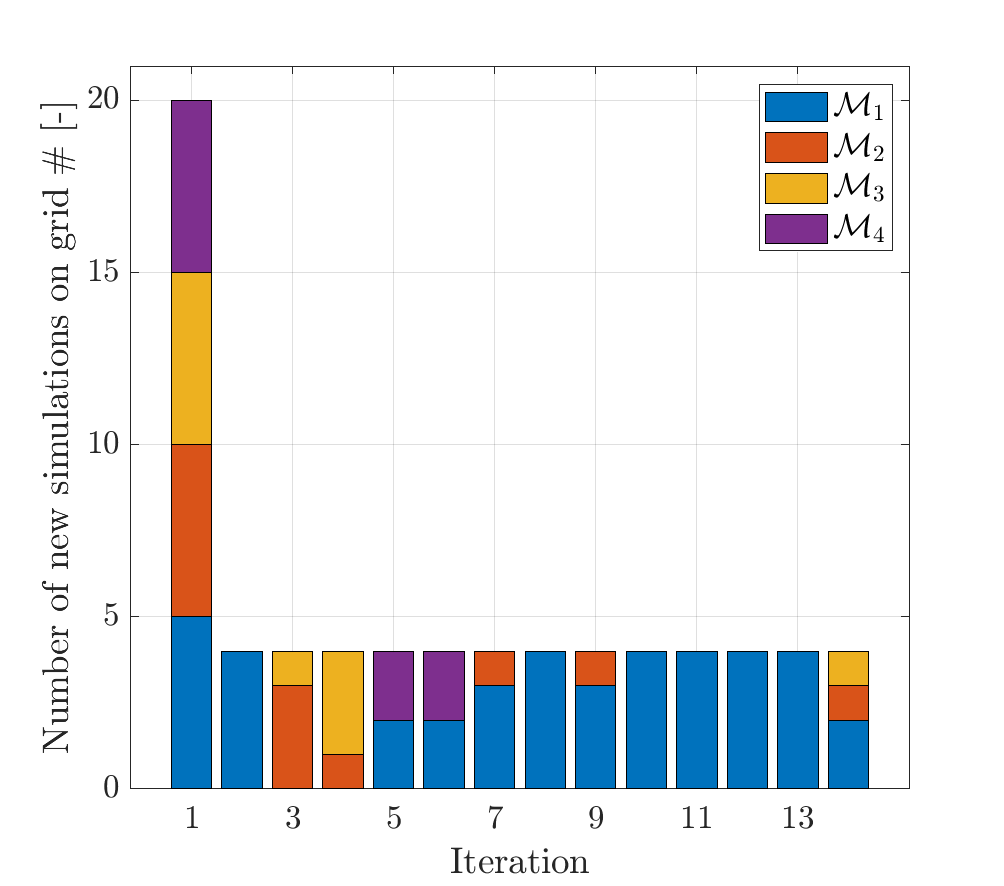}}
	\caption{Comparison of MISC and SRBF results: number of simulations required at each iteration for each grid. }
	\label{fig:beta_points_comparison}
\end{figure}

% RISULTATI SRBF
\begin{figure}[!t]
	\centering
	\subfigure[]{\includegraphics[width=1\textwidth]{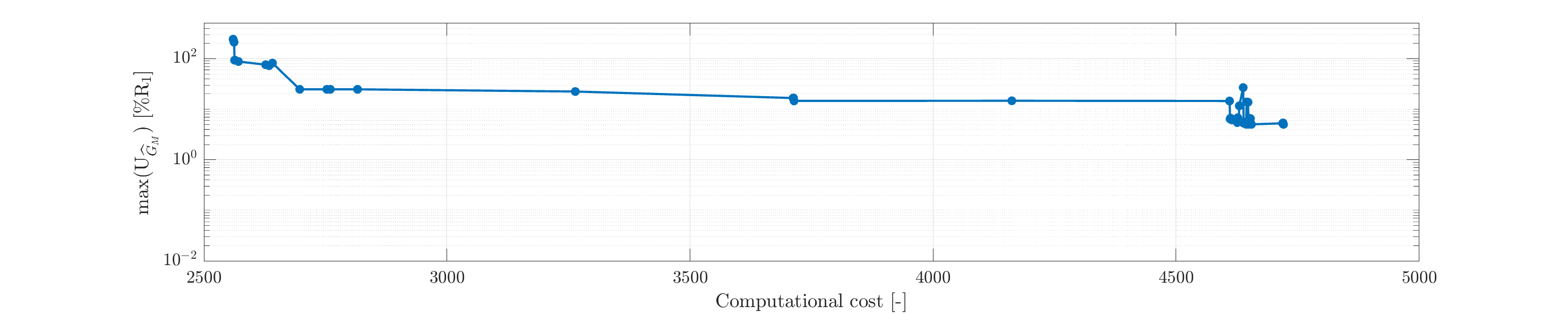}}
	\subfigure[]{\includegraphics[width=1\textwidth]{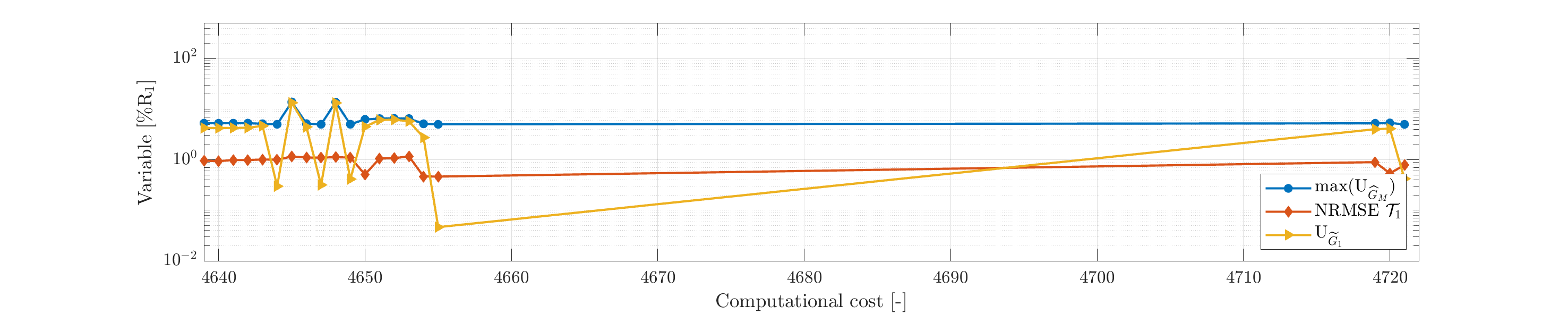}}
	\caption{SRBF results: (a) convergence of the maximum prediction uncertainty and (b) its comparison with the lowest-fidelity prediction uncertainty and the normalized (with the function range) RMSE of $\mathcal{T}_1$ versus the computational cost from iteration $8$ to final.}
	\label{fig:points_RBF}
\end{figure}

\begin{figure}[!t]
	\centering
	\subfigure[MISC, iteration 14]{\includegraphics[width=0.45\textwidth]{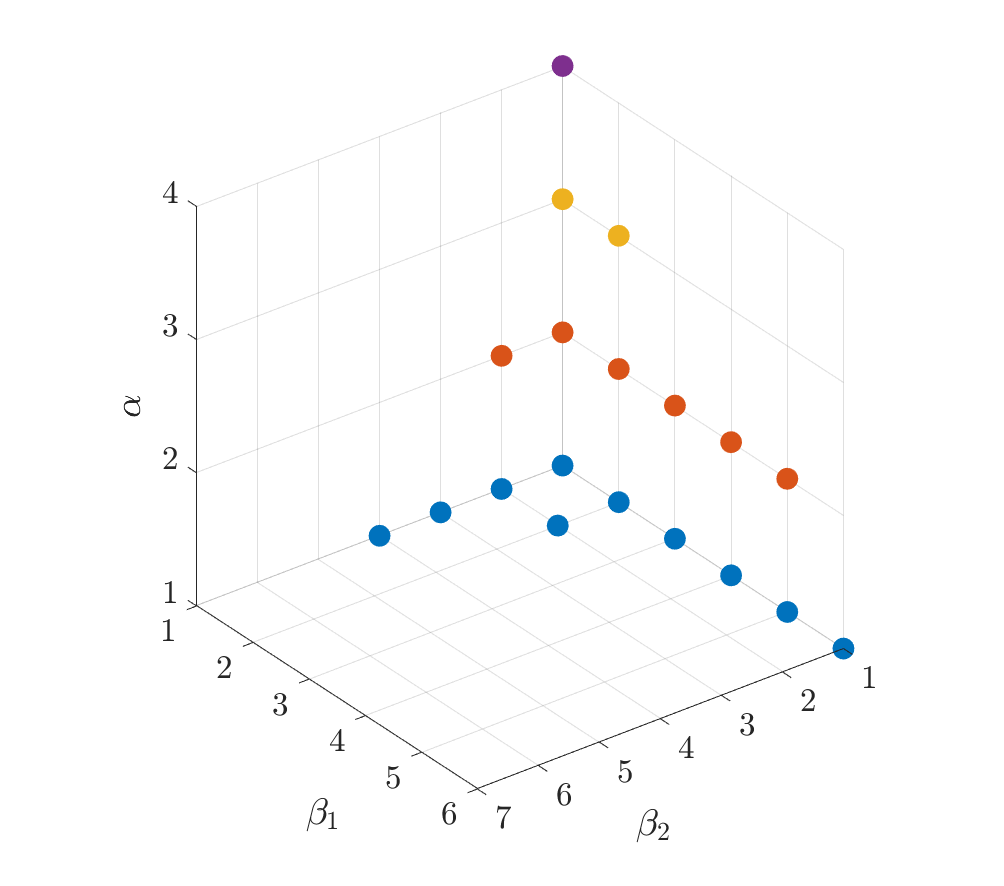}}
	\subfigure[MISC, final iteration]{\includegraphics[width=0.45\textwidth]{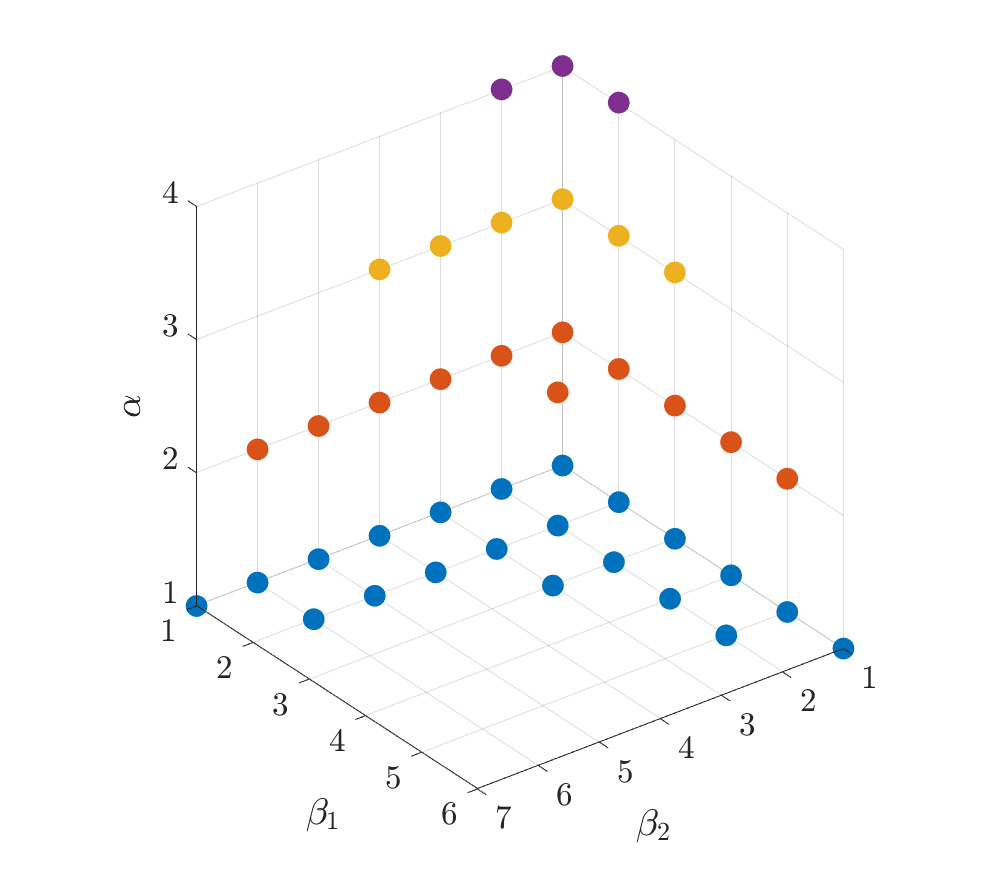}}
	\caption{MISC results: explored multi-indices. The indices $\beta_1$ and $\beta_2$ determine how many quadrature points are selected according to the level function \eqref{eq:lev_fun} for the advancement speed and the draught, respectively.}
	\label{fig:indices_MISC}
\end{figure}

\begin{figure}[!t]
	\centering
	\includegraphics[width=1\textwidth]{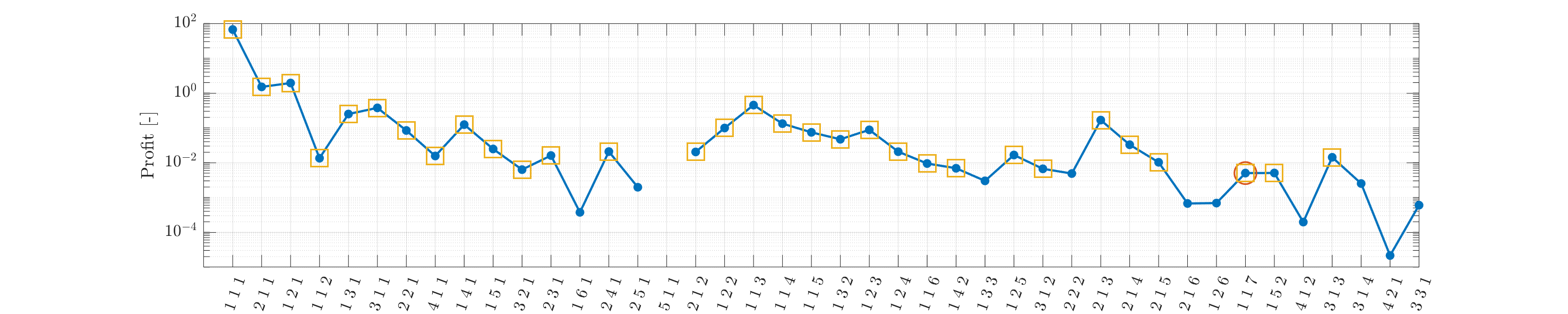}
	\caption{MISC results: profit of the explored multi-indices at the final iteration. The multi-indices circled in yellow have profit higher than the one marked in red and have been added to the index set $\bm{I}$ (see Algorithm \ref{algo:misc_implementation}).}
	\label{fig:profits_MISC}
\end{figure}

Fig. \ref{fig:EV_comparison} shows the convergence of the expected value and standard deviation of $R_T$ versus the computational cost. MISC and SRBF converge towards similar estimates. MISC achieves good estimate already with low computational cost, whereas at later iterations, corresponding to computational cost above 1000, its behavior worsens; conversely SRBF gives a more stable estimate. Similar conclusions can be drawn by looking at the plot of the values of the standard deviation. Complementary information are given in Tab. \ref{table:comparison}, summarizing the expected value and the standard deviation achieved by the two methods at the final iteration and at one characteristic iteration. In the case of MISC iteration 14 (associated to a computational cost equals to 900) is selected, one of the iterations giving a good estimate of the expected value but before the ``disruption'' happening around cost approximately 1000. Such iteration is indicative of the behavior of MISC and the successive discussions are referred to it.
{Concerning SRBF, the} convergence of the expected value and the standard deviation is oscillatory (zoom not shown for brevity). The oscillation becomes more evident after iteration $8$ with computational cost of $4638$, {which is then selected as characteristic iteration for the SRBF and shown in Tab. \ref{table:comparison}). Actually,} after this iteration the SRBF moves from interpolation to {least-squares} approximation for the lowest fidelity only, since its number of training points reaches the threshold value of $5^N$. This produces oscillations in the last part of the SRBF convergence. Nevertheless, variations are small and oscillations limited.

Fig. \ref{fig:pdf_comparison} shows the $R_T$ histogram and {the estimate of the PDF of $R_T$ obtained using a kernel density estimation method}.
{Crucially, the values of $R_T$ used build histograms and the approximate PDF can be obtained cheaply, by evaluating the MISC and SRBF surrogate models over a suitable number of Monte Carlo samples from the parameter space.} 
The PDF obtained with SRBF has a peak around in the interval 20-40. A similar peak is retrieved also by MISC. However, in the latter case secondary peaks are observed for higher values of resistance at the final iteration. SRBF with exact interpolation (see Fig. \ref{fig:pdf_comparison}c) leads to a peak value significantly higher than MISC and SRBF with least-squares regression (see Fig. \ref{fig:pdf_comparison}a, b, and d).
{The spurious peaks in the PDF obtained with MISC can be explained by looking at the MISC response surface for $R_T$, which we discuss next.}

{The response surfaces obtained by the two methods} are shown and compared in Fig. \ref{fig:response_surface_comparison}. Once again, it can be observed that the behavior of MISC deteriorates at later iterations. The response surface is very irregular and this can be attributed to the numerical noise which affects the CFD simulations {(due to the fact that the solver implements an iterative method
  which stops as soon as a prescribed tolerance is met)}, in particular the one on the coarsest grid.
As the overall idea of MISC is to solve most PDEs on the less expensive grids and building the surrogate model with lagrangian (hence exact) interpolation, 
the presence of numerical noise turns out to be problematic for this method. This issue could be mitigated by doing e.g. a least-squares regression on the data, to smooth the numerical oscillations. A more thorough investigation of this aspect is out of scope here and will be considered in a future work (see also discussion in Sect. \ref{sect:conclusions} - Conclusions).
% One could also take advantage of the soft information (monotonicity, multimodality, ...) on the physical nature of the problem. In this particular case, the resistance is expected to be monotonically increasing with respect to advancement speed and draught. Such property could be preserved by employing, e.g., least-squares regressions with appropriate polynomial degree and/or monotonic smoothing \cite{dykstra1982}. Investigating this aspect is out of scope here and it will be considered in a future work.
On the contrary, SRBF suffers the numerical noise of the CFD outputs {only as long as} exact interpolation is imposed (see Fig. \ref{fig:response_surface_comparison}c) leading to a significant deformation of the response surface close to the lower bound of both the uncertain parameters,
whereas {for later iterations (i.e., larger computational costs)}
the use of regression improves the quality of the surrogates by smoothing the response surface
and filtering out the numerical noise (see Fig. \ref{fig:response_surface_comparison}d). 

Fig. \ref{fig:points_comparison} displays the points in the uncertain parameter space selected by MISC and SRBF. The two plots on the top line show the points employed by MISC at iteration 14 and at the final iteration. At the beginning the MISC algorithm explores more {the variability of the advancement speed}, suggesting a stronger dependence of the quantity of interest on {this parameter} rather than on the draught.
Points exploring {the variability of the draught} are added only at later iterations (this subject is elaborated also later on).
The sampling performed by SRBF is shown in Figs. \ref{fig:points_comparison}c and Figs. \ref{fig:points_comparison}d {at iteration 8} and {at the} final iteration, respectively. Until iteration $8$ (exact interpolation only), the adaptive sampling method explores the domain extrema and starts to cluster samples in two zones among $U=[7,8.8]$ m/s, since the numerical noise negatively affects the prediction uncertainty of the interpolating SRBF. Switching to the least-squares approximation and filtering out the noise prevents an excessive clusterization of the samples at the successive iterations. Indeed, the SRBF samples are fairly spread over the domain. 

Fig. \ref{fig:beta_points_comparison} shows the amount of simulations on the different spatial grids required by both methods at each iteration. Specifically, on the left the total number of simulations is displayed, whereas on the right only the number of new simulations asked at each iteration is shown. MISC required most of the simulations on the coarse grid; only five simulations on the finest grid are required (one at iteration 5, two at iteration 28, and two at iteration 30). 
At iteration 12, 24, and 25 no new points are added, which can be explained as follows. At iteration 12 the MISC algorithm requires simulations for spatial index $\alpha=5$: however, having considered only four spatial refinements level, this grid was not available. The profit of this multi-index is then set to zero and the algorithm restarts to explore from the multi-index with the next highest profit. 
Instead, at iteration 24 and 25 there is no new index $\jj$ to be added to the multi-index set because none of the indices $\jj \in \text{Mar}(\bm{I})$ would lead to $\bm{I} \cup \{\jj\}$ downward closed, see Algorithm \ref{algo:misc_implementation}. As for SRBF, it {spent} about 50\% of the final computational cost at the first iteration, then requiring simulations on the finest grids only at iterations 5 and 6. In all the other iterations mainly low-fidelity simulations are performed. This sampling behavior is due to the high values of prediction uncertainty that are found in the corners of the variables domain, because the topology of the initial training leads to extrapolation in those zones. 
%In that regions the highest prediction uncertainty is estimated and the adaptive sampling procedure have required all the fidelities, before moving to explore other regions.
Such corner regions are those with the highest estimated prediction uncertainty,
and the adaptive sampling procedure requires all the fidelities before moving
to explore other regions.

Fig. \ref{fig:points_RBF}a shows the maximum prediction uncertainty versus the computational cost. After the first four iterations the value of the maximum prediction uncertainty decreases from 120\% of the high-fidelity function range (R$_1$) to about 25\%. During the adaptive sampling the value of the maximum prediction uncertainty continues to decrease slower than in the first iterations, achieving a final value of about 5\%. The convergence shows some oscillations around a computational cost equal to 4600, due to the change from exact interpolation to least-squares approximation of the training set. 
Fig. \ref{fig:points_RBF}b compares the convergence of the maximum prediction uncertainty with its low-fidelity component and the normalized (with the function range) RMSE of the lowest-fidelity training set from iteration $8$. It is worth noting that the lowest-fidelity prediction uncertainty is almost equal to the maximum multi-fidelity prediction uncertainty until a computational cost of about 4655, thus representing its main component. %After lowest-fidelity samples are added, the lowest-fidelity components decrease significantly and a $\mathcal{M}_3$ simulation is requested. Therefore, the $U_{\hat{\varepsilon}_3}$ component of the prediction uncertainty is reduced and the low-fidelity metamodel prediction uncertainty component is again the greatest. Finally, the NRMSE of the fourth fidelity training set is almost constant and its value is significantly smaller than the maximum prediction uncertainty.

% RISULTATI MISC

Finally, coming back to the analysis of the MISC results, in Fig. \ref{fig:indices_MISC} the multi-indices $\left[\alpha, \beta_1,\beta_2\right]$ selected by the MISC method until iteration 14 and at the final iteration are shown. It can be seen that the multi-index set at iteration 14 stretches more in direction ``advancement speed''
{($\beta_1$)} than in direction ``draught'' {($\beta_2$)}:
this further confirms the statement already made that the parameter advancement speed is deemed more influential by MISC, i.e. the multi-indices with larger profits are predominantly found in direction advancement speed. More quantitative information on the profits can be found in Fig. \ref{fig:profits_MISC}, where the profit of every explored multi-index is displayed. Note that only the indices with the highest profits (marked in yellow) are added to the multi-index set $\bm{I}$. The remaining indices are explored (i.e. added to the multi-index set $\bm{G}$, therefore contribute to the MISC approximation) but non yet selected (i.e. added to the multi-index set $\bm{I}$). As already explained earlier on, there is no profit associated to the multi-index $\left[5,1,1\right]$ since the fifth spatial refinements level is not available.

\section{Conclusions and Future Work}\label{sect:conclusions}
MISC and SRBF have been applied and compared for the UQ of a RoPax ferry advancing in calm water and subject to two operational uncertainties, namely the ship speed and draught (linked to the payload). The estimation of expected value, standard deviation, and probability density function of the (model-scale) resistance obtained by RANS computations has been discussed. 
Both MISC and SRBF used as multi-fidelity levels the multi-grid computations provided by the RANS solver.
More specifically, four grids (obtained as isotropic coarsening of an initial grid) were available, and have been used.
This implies that the multi-indices considered in MISC have just one component controlling
the spatial discretization (instead of three), so that in total the dimension of the multi-indices is three (one component for the spatial discretization
plus two components for the parametric discretization).

  Overall, MISC proved to be quite effective in delivering a reasonable approximation of the statistical indicators (expected value, standard deviation, probability
  density function) of the quantity of interest with a limited computational cost, but on the other hand turned out to be
  very sensitive to the numerical noise that affects the RANS solver, especially with the coarsest grid; this becomes clearly visible when the computational budget gets
  larger, as more runs of the RANS solver are considered. Strategies to mitigate this effect
  by smoothing the data coming from the RANS solver (e.g. by least-squares regression) are under investigation.
  Such smoothing could also take into account the soft information (monotonicity, multimodality, etc.)
  available on the physical nature of the problem. In this particular case, the resistance is expected to
  be monotone increasing with respect to advancement speed and draught: such property could be preserved
  by employing, e.g., least-squares regressions with appropriate polynomial degrees and/or monotonic smoothing, see e.g. \cite{dykstra1982}.
  Another practical problem is caused by the non-monotonic behavior of the profits,
  where some indices with low-profits shade useful neighbors, thus delay the convergence of MISC. 
  More robust strategies to explore the set of multi-indices, that blend the profit-based selection of indices with other criteria
  are also subject of current work; see e.g. \cite{chkifa:adaptive-interp,gerstner.griebel:adaptive}, where this problem
  was discussed in the context of adaptive sparse-grids quadrature/interpolation.

%SRBF
%
% Conclusioni:
% - Il capionamento MF come desiderato predilige le valutazione a fedelta' piu' bassa
% - Costo iniziale del training e' abbastanza elevato, ciononostante la convergenza e' abbastanza robusta 
% - L'approccio interpolante poi regressivo sembra riuscire a trattare in modo soddisfacente il rumore
%
% Lavoro futuro:
% Campionamento iniziale
% Investigare sul criterio di attivazione della regressione attivazione
% GP
% Sampling w/ w/o considering noise (questo non lo metto, troppo specifico a GP)
%
The SRBF training also used mainly lowest-fidelity (and therefore lowest computational cost) RANS computations, as desired. Nevertheless, the initial training (first iteration) had a quite high computational cost, due to the use of RANS computations from all fidelity levels. Within the current computational effort, the convergence of expected value, standard deviation, and probability density function of the ship resistance is reasonably smooth, with quite small variations along the iteration. Moreover, the adaptive use of regressive SRBF has a beneficial effect on the robustness of the method to noisy data. 
%Aggiungere griglia rada = lot of noise
Ongoing and future work focuses on the definition of a more efficient scheme for the initial training, as well as a data-adaptive criterion for the activation of the regressive model. Comparison with other methods, such as multi-fidelity Gaussian processes, will be also investigated.

% Conclusioni generali
% Per poche valutazioni meglio MISC
% per piu' valutazioni equivalenti con una leggere preferenza per RBF docuta alla capacita' di trattamento del rumore
%
% Future work generale
% Necessita' valori d benchmark per stima degli errori su ev, std, pdf
% Applicazione a casi piu' complessi con maggiore dimensionalita'
%
MISC and SRBF give close values of expected value and standard deviation. The probability density functions are reasonably close, although the MISC surrogate model is affected by the noise in the data and its effects are visible in the resulting density function. Overall, the results suggest that MISC could be preferred when only a limited number of RANS computations is available. For larger data sets both MISC and SRBF represent a valid option, with a slight preference in favor of the current implementation of SRBF, due to its robustness to noise. Future research will address more complex test cases (larger number of uncertain parameters and more realistic conditions, such as regular/irregular waves) possibly validating the results against benchmark values.

\section*{Acknowledgments}
CNR-INM is grateful to Dr. Woei-Min Lin, Dr. Elena McCarthy, and Dr. Salahuddin Ahmed of the Office of Naval Research and Office of Naval Research Global, for their support through NICOP grant N62909-18-1-2033. Dr. Riccardo Pellegrini is partially
supported through CNR-INM project OPTIMAE. The HOLISHIP project (HOLIstic optimisation of SHIP design and operation for life cycle, www.holiship.eu) is also acknowledged, funded by the European Union’s Horizon 2020 research and innovation program under grant agreement N. 689074. 
Lorenzo Tamellini and Chiara Piazzola have been supported by the PRIN 2017 project 201752HKH8 ``Numerical Analysis for Full and Reduced Order Methods for the efficient and accurate solution of complex systems governed by Partial Differential Equations (NA-FROM-PDEs)''.
Lorenzo Tamellini also acknowledges the support of GNCS-INdAM (Gruppo Nazionale Calcolo Scientifico - Istituto Nazionale di Alta Matematica).

\bibliography{biblio}

\begin{thebibliography}{45}
\newcommand{\enquote}[1]{``#1''}
\providecommand{\natexlab}[1]{#1}
\providecommand{\url}[1]{\texttt{#1}}
\providecommand{\urlprefix}{URL }
\expandafter\ifx\csname urlstyle\endcsname\relax
  \providecommand{\doi}[1]{\discretionary{}{}{}https://doi.org/#1}\else
  \providecommand{\doi}[1]{\discretionary{}{}{}\urlstyle{rm}\url{https://doi.org/#1}}\fi

\bibitem[{He et~al.(2013)He, Diez, Zou, Campana, and Stern}]{he2013urans}
He, W., Diez, M., Zou, Z., Campana, E.~F., and Stern, F., \enquote{URANS study
  of Delft catamaran total/added resistance, motions and slamming loads in head
  sea including irregular wave and uncertainty quantification for variable
  regular wave and geometry,} \emph{Ocean Engineering}, Vol.~74, 2013, pp.
  189--217.

\bibitem[{Diez et~al.(2018{\natexlab{a}})Diez, Broglia, Durante, Olivieri,
  Campana, and Stern}]{diez2018statistical}
Diez, M., Broglia, R., Durante, D., Olivieri, A., Campana, E.~F., and Stern,
  F., \enquote{Statistical Assessment and Validation of Experimental and
  Computational Ship Response in Irregular Waves,} \emph{Journal of
  Verification, Validation and Uncertainty Quantification}, Vol.~3, No.~2,
  2018{\natexlab{a}}.

\bibitem[{Durante et~al.(2020)Durante, Broglia, Diez, Olivieri, Campana, and
  Stern}]{durante2020accurate}
Durante, D., Broglia, R., Diez, M., Olivieri, A., Campana, E., and Stern, F.,
  \enquote{Accurate experimental benchmark study of a catamaran in regular and
  irregular head waves including uncertainty quantification,} \emph{Ocean
  Engineering}, Vol. 195, 2020, p. 106685.

\bibitem[{Diez et~al.(2018{\natexlab{b}})Diez, Campana, and
  Stern}]{diez2018stochastic}
Diez, M., Campana, E.~F., and Stern, F., \enquote{Stochastic optimization
  methods for ship resistance and operational efficiency via CFD,}
  \emph{Structural and Multidisciplinary Optimization}, Vol.~57, No.~2,
  2018{\natexlab{b}}, pp. 735--758.

\bibitem[{Serani et~al.(2019{\natexlab{a}})Serani, Diez, Wackers, Visonneau,
  and Stern}]{serani2019stochastic}
Serani, A., Diez, M., Wackers, J., Visonneau, M., and Stern, F.,
  \enquote{Stochastic shape optimization via design-space augmented
  dimensionality reduction and rans computations,} \emph{57th AIAA Aerospace
  Sciences Meeting, {SciTech} 2019}, 2019{\natexlab{a}}, p. 2218.

\bibitem[{Quagliarella et~al.(2019)Quagliarella, Serani, Diez, Pisaroni,
  Leyland, Montagliani, Iemma, Gaul, Shin, Wunsch, Hirsch, Choi, and
  Stern}]{quagliarella2019benchmarking}
Quagliarella, D., Serani, A., Diez, M., Pisaroni, M., Leyland, P., Montagliani,
  L., Iemma, U., Gaul, N.~J., Shin, J., Wunsch, D., Hirsch, C., Choi, K., and
  Stern, F., \enquote{Benchmarking Uncertainty Quantification Methods Using the
  NACA 2412 Airfoil with Geometrical and Operational Uncertainties,} \emph{57th
  AIAA Aerospace Sciences Meeting, {SciTech} 2019}, 2019, p. 3555.

\bibitem[{Beran et~al.(2020)Beran, Bryson, Thelen, Diez, and
  Serani}]{beran2020}
Beran, P.~S., Bryson, D.~E., Thelen, A.~S., Diez, M., and Serani, A.,
  \enquote{Comparison of Multi-Fidelity Approaches for Military Vehicle
  Design,} \emph{21th AIAA/ISSMO Multidisciplinary Analysis and Optimization
  Conference ({MA\&O}), {AVIATION} 2020}, Virtual Event, June 15-19, 2020.

\bibitem[{Giles(2008)}]{giles:MLMC}
Giles, M.~B., \enquote{Multilevel Monte Carlo Path Simulation,}
  \emph{Operations Research}, Vol.~56, No.~3, 2008, pp. 607--617.

\bibitem[{Cliffe et~al.(2011)Cliffe, Giles, Scheichl, and
  Teckentrup}]{scheichl.giles:MLMC}
Cliffe, K., Giles, M., Scheichl, R., and Teckentrup, A., \enquote{Multilevel
  Monte Carlo methods and applications to elliptic PDEs with random
  coefficients,} \emph{Computing and Visualization in Science}, Vol.~14, No.~1,
  2011, pp. 3--15.

\bibitem[{Haji-Ali et~al.(2015)Haji-Ali, Nobile, and
  Tempone}]{hajiali.eal:MultiIndexMC}
Haji-Ali, A.-L., Nobile, F., and Tempone, R., \enquote{{M}ulti-{i}ndex {M}onte
  {C}arlo: when sparsity meets sampling,} \emph{Numerische Mathematik}, 2015,
  pp. 1--40.

\bibitem[{Kuo et~al.(2015)Kuo, Schwab, and Sloan}]{kss12}
Kuo, F.~Y., Schwab, C., and Sloan, I., \enquote{Multi-level {Q}uasi-{M}onte
  {C}arlo {F}inite {E}lement {M}ethods for a {C}lass of {E}lliptic {PDE}s with
  {R}andom {C}oefficients,} \emph{Foundations of Computational Mathematics},
  Vol.~15, No.~2, 2015, pp. 411--449.

\bibitem[{Teckentrup et~al.(2015)Teckentrup, Jantsch, Webster, and
  Gunzburger}]{teckentrup.etal:MLSC}
Teckentrup, A.~L., Jantsch, P., Webster, C.~G., and Gunzburger, M., \enquote{A
  {M}ultilevel {S}tochastic {C}ollocation {M}ethod for {P}artial {D}ifferential
  {E}quations with {R}andom {I}nput {D}ata,} \emph{SIAM/ASA Journal on
  Uncertainty Quantification}, Vol.~3, No.~1, 2015, pp. 1046--1074.

\bibitem[{Beck et~al.(2019)Beck, Tamellini, and Tempone}]{beck.eal:MISC-IGA}
Beck, J., Tamellini, L., and Tempone, R., \enquote{{IGA-based Multi-Index
  Stochastic Collocation for random PDEs on arbitrary domains},}
  \emph{{Computer Methods in Applied Mechanics and Engineering}}, Vol. 351,
  2019, pp. 330--350.

\bibitem[{Jakeman et~al.(2020)Jakeman, Eldred, Geraci, and
  Gorodetsky}]{jakeman2019adaptive}
Jakeman, J.~D., Eldred, M., Geraci, G., and Gorodetsky, A., \enquote{Adaptive
  Multi-index Collocation for Uncertainty Quantification and Sensitivity
  Analysis,} \emph{International Journal for Numerical Methods in Engineering},
  Vol. 121, No.~6, 2020, pp. 1314--1343.

\bibitem[{Haji-Ali et~al.(2016{\natexlab{a}})Haji-Ali, Nobile, Tamellini, and
  Tempone}]{hajiali.eal:MISC1}
Haji-Ali, A., Nobile, F., Tamellini, L., and Tempone, R., \enquote{Multi-Index
  Stochastic Collocation for random {PDEs},} \emph{Computer Methods in Applied
  Mechanics and Engineering}, Vol. 306, 2016{\natexlab{a}}, pp. 95--122.
\newblock \doi{http://dx.doi.org/10.1016/j.cma.2016.03.029}.

\bibitem[{Haji-Ali et~al.(2016{\natexlab{b}})Haji-Ali, Nobile, Tamellini, and
  Tempone}]{hajiali.eal:MISC2}
Haji-Ali, A.-L., Nobile, F., Tamellini, L., and Tempone, R.,
  \enquote{Multi-index {S}tochastic {C}ollocation convergence rates for random
  {PDE}s with parametric regularity,} \emph{Foundations of {C}omputational
  {M}athematics}, Vol.~16, No.~6, 2016{\natexlab{b}}, pp. 1555--1605.
\newblock \doi{10.1007/s10208-016-9327-7}.

\bibitem[{{Haji-Ali, Abdul-Lateef} et~al.(2020){Haji-Ali, Abdul-Lateef},
  {Nobile, Fabio}, {Tempone, Ra\'ul}, and {Wolfers,
  S\"oren}}]{hajiali2017multilevel}
{Haji-Ali, Abdul-Lateef}, {Nobile, Fabio}, {Tempone, Ra\'ul}, and {Wolfers,
  S\"oren}, \enquote{Multilevel weighted least squares polynomial
  approximation,} \emph{ESAIM: M2AN}, Vol.~54, No.~2, 2020, pp. 649--677.

\bibitem[{Peherstorfer et~al.(2018)Peherstorfer, Willcox, and
  Gunzburger}]{peherstorfer:MFsurvey}
Peherstorfer, B., Willcox, K., and Gunzburger, M., \enquote{Survey of
  Multifidelity Methods in Uncertainty Propagation, Inference, and
  Optimization,} \emph{SIAM Review}, Vol.~60, No.~3, 2018, pp. 550--591.

\bibitem[{Gorodetsky et~al.(2020)Gorodetsky, Geraci, Eldred, and
  Jakeman}]{Gorodetsky2020b}
Gorodetsky, A., Geraci, G., Eldred, M.~S., and Jakeman, J., \enquote{A
  generalized approximate control variate framework for multifidelity
  uncertainty quantification,} \emph{Journal of Computational Physics}, Vol.
  408, 2020, p. 109257.

\bibitem[{Pisaroni et~al.(2017{\natexlab{a}})Pisaroni, Nobile, and
  Leyland}]{Pisaroni:MLMC-opt}
Pisaroni, M., Nobile, F., and Leyland, P., \enquote{A Multilevel Monte Carlo
  Evolutionary Algorithm for Robust Aerodynamic Shape Design,} \emph{18th
  AIAA/ISSMO Multidisciplinary Analysis and Optimization Conference. Denver,
  Colorado}, 2017{\natexlab{a}}.

\bibitem[{Pisaroni et~al.(2017{\natexlab{b}})Pisaroni, Nobile, and
  Leyland}]{pisaroni:CMLMC}
Pisaroni, M., Nobile, F., and Leyland, P., \enquote{A Continuation Multi Level
  Monte Carlo (C-MLMC) method for uncertainty quantification in compressible
  inviscid aerodynamics,} \emph{Computer Methods in Applied Mechanics and
  Engineering}, Vol. 326, 2017{\natexlab{b}}, pp. 20--50.

\bibitem[{Geraci et~al.(2019)Geraci, Eldred, Gorodetsky, and
  Jakeman}]{geraci2019}
Geraci, G., Eldred, M.~S., Gorodetsky, A., and Jakeman, J., \enquote{Recent
  advancements in Multilevel-Multifidelity techniques for forward UQ in the
  DARPA Sequoia project,} \emph{57th AIAA Aerospace Sciences Meeting, {SciTech}
  2019}, 2019, p. 0722.

\bibitem[{Han and G{\"o}rtz(2012)}]{han2012-AIAA}
Han, Z.-H., and G{\"o}rtz, S., \enquote{Hierarchical kriging model for
  variable-fidelity surrogate modeling,} \emph{AIAA journal}, Vol.~50, No.~9,
  2012, pp. 1885--1896.

\bibitem[{Baar et~al.(2015)Baar, Roberts, Dwight, and Mallol}]{debaar2015-CF}
Baar, J.~d., Roberts, S., Dwight, R., and Mallol, B., \enquote{Uncertainty
  quantification for a sailing yacht hull, using multi-fidelity kriging,}
  \emph{Computers \& Fluids}, Vol. 123, 2015, pp. 185--201.

\bibitem[{Wackers et~al.(2020{\natexlab{a}})Wackers, Visonneau, Pellegrini,
  Ficini, Serani, and Diez}]{wackers2020-AIAA}
Wackers, J., Visonneau, M., Pellegrini, R., Ficini, S., Serani, A., and Diez,
  M., \enquote{Adaptive N-Fidelity Metamodels for Noisy {CFD} Data,} \emph{21th
  AIAA/ISSMO Multidisciplinary Analysis and Optimization Conference ({MA\&O}),
  {AVIATION} 2020}, Virtual Event, June 15-19, 2020{\natexlab{a}}.

\bibitem[{Serani et~al.(2019{\natexlab{b}})Serani, Pellegrini, Wackers,
  Jeanson, Queutey, Visonneau, and Diez}]{serani2019-IJCFD}
Serani, A., Pellegrini, R., Wackers, J., Jeanson, C.-J., Queutey, P.,
  Visonneau, M., and Diez, M., \enquote{Adaptive multi-fidelity sampling for
  {CFD}-based optimization via radial basis functions metamodel,}
  \emph{International Journal of Computational Fluid Dynamics}, Vol.~33, No.
  6-7, 2019{\natexlab{b}}, pp. 237--255.

\bibitem[{Han et~al.(2013)Han, G{\"o}rtz, and Zimmermann}]{han2013-AST}
Han, Z.-H., G{\"o}rtz, S., and Zimmermann, R., \enquote{Improving
  variable-fidelity surrogate modeling via gradient-enhanced kriging and a
  generalized hybrid bridge function,} \emph{Aerospace Science and Technology},
  Vol.~25, No.~1, 2013, pp. 177--189.

\bibitem[{Di~Mascio et~al.(2007)Di~Mascio, Broglia, and Muscari}]{dimascio2007}
Di~Mascio, A., Broglia, R., and Muscari, R., \enquote{On the application of the
  single-phase level set method to naval hydrodynamic flows,} \emph{Computers
  \& fluids}, Vol.~36, No.~5, 2007, pp. 868--886.

\bibitem[{Di~Mascio et~al.(2009)Di~Mascio, Broglia, and Muscari}]{dimascio2009}
Di~Mascio, A., Broglia, R., and Muscari, R., \enquote{Prediction of
  hydrodynamic coefficients of ship hulls by high-order Godunov-type methods,}
  \emph{Journal of Marine Science and Technology}, Vol.~14, No.~1, 2009, pp.
  19--29.

\bibitem[{Broglia and Durante(2018)}]{broglia2018}
Broglia, R., and Durante, D., \enquote{Accurate prediction of complex free
  surface flow around a high speed craft using a single-phase level set
  method,} \emph{Computational Mechanics}, Vol.~62, No.~3, 2018, pp. 421--437.

\bibitem[{Trefethen(2008)}]{trefethen:comparison}
Trefethen, L.~N., \enquote{Is {G}auss quadrature better than
  {C}lenshaw-{C}urtis?} \emph{SIAM Rev.}, Vol.~50, No.~1, 2008, pp. 67--87.

\bibitem[{Gerstner and Griebel(2003)}]{gerstner.griebel:adaptive}
Gerstner, T., and Griebel, M., \enquote{Dimension-adaptive tensor-product
  quadrature,} \emph{Computing}, Vol.~71, No.~1, 2003, pp. 65--87.
\newblock \doi{10.1007/s00607-003-0015-5},
  \urlprefix\url{http://dx.doi.org/10.1007/s00607-003-0015-5}.

\bibitem[{Nobile et~al.(2016)Nobile, Tamellini, Tesei, and
  Tempone}]{nobile.eal:adaptive-lognormal}
Nobile, F., Tamellini, L., Tesei, F., and Tempone, R., \enquote{An adaptive
  sparse grid algorithm for elliptic {PDE}s with lognormal diffusion
  coefficient,} \emph{Sparse Grids and Applications -- Stuttgart 2014}, Lecture
  Notes in Computational Science and Engineering, Vol. 109, edited by J.~Garcke
  and D.~Pfl\"uger, Springer International Publishing Switzerland, 2016, pp.
  191--220.

\bibitem[{Klimke(2006)}]{klimke:thesis}
Klimke, A., \enquote{Uncertainty modeling using fuzzy arithmetic and sparse
  grids,} Ph.D. thesis, Universit\"at Stuttgart, Shaker Verlag, Aachen, 2006.

\bibitem[{Chkifa et~al.(2014)Chkifa, Cohen, and
  Schwab}]{chkifa:adaptive-interp}
Chkifa, A., Cohen, A., and Schwab, C., \enquote{High-Dimensional Adaptive
  Sparse Polynomial Interpolation and Applications to Parametric PDEs,}
  \emph{Foundations of Computational Mathematics}, Vol.~14, No.~4, 2014, pp.
  601--633.
\newblock \doi{10.1007/s10208-013-9154-z}.

\bibitem[{Guignard and Nobile(2018)}]{Guignard:a-post}
Guignard, D., and Nobile, F., \enquote{A Posteriori Error Estimation for the
  Stochastic Collocation Finite Element Method,} \emph{SIAM Journal on
  Numerical Analysis}, Vol.~56, No.~5, 2018, pp. 3121--3143.

\bibitem[{Gutmann(2001)}]{gutmann2001}
Gutmann, H.~M., \enquote{A radial basis function method for global
  optimization,} \emph{Journal of global optimization}, Vol.~19, No.~3, 2001,
  pp. 201--227.

\bibitem[{Forrester and Keane(2009)}]{forrester2009}
Forrester, A. I.~J., and Keane, A.~J., \enquote{Recent advances in
  surrogate-based optimization,} \emph{Progress in aerospace sciences},
  Vol.~45, No. 1-3, 2009, pp. 50--79.

\bibitem[{Volpi et~al.(2015)Volpi, Diez, Gaul, Song, Iemma, Choi, Campana, and
  Stern}]{volpi2015-SMO}
Volpi, S., Diez, M., Gaul, N.~J., Song, H., Iemma, U., Choi, K.~K., Campana,
  E.~F., and Stern, F., \enquote{Development and validation of a dynamic
  metamodel based on stochastic radial basis functions and uncertainty
  quantification,} \emph{Structural and Multidisciplinary Optimization},
  Vol.~51, No.~2, 2015, pp. 347--368.

\bibitem[{Lloyd(1982)}]{lloyd1982-IEEE}
Lloyd, S., \enquote{Least squares quantization in PCM,} \emph{IEEE transactions
  on information theory}, Vol.~28, No.~2, 1982, pp. 129--137.

\bibitem[{Li et~al.(2017)Li, Gao, Gu, Gong, Jing, and Su}]{li2017-SMO}
Li, X., Gao, W., Gu, L., Gong, C., Jing, Z., and Su, H., \enquote{A cooperative
  radial basis function method for variable-fidelity surrogate modeling,}
  \emph{Structural and Multidisciplinary Optimization}, Vol.~56, No.~5, 2017,
  pp. 1077--1092.

\bibitem[{Serani et~al.(2019{\natexlab{c}})Serani, Pellegrini, Broglia,
  Wackers, Visonneau, and Diez}]{serani2019-MARINE}
Serani, A., Pellegrini, R., Broglia, R., Wackers, J., Visonneau, M., and Diez,
  M., \enquote{An Adaptive N-Fidelity Metamodel for Design and
  Operational-Uncertainty Space Exploration of Complex Industrail Problems,}
  \emph{Proceedings of the 8th International Conference on Computational
  Methods in Marine Engineering (Marine 2019)}, 2019{\natexlab{c}}, pp.
  177--188.

\bibitem[{Wackers et~al.(2020{\natexlab{b}})Wackers, Visonneau, Serani,
  Pellegrini, Broglia, and Diez}]{wackers2020-SNH}
Wackers, J., Visonneau, M., Serani, A., Pellegrini, R., Broglia, R., and Diez,
  M., \enquote{Multi-Fidelity Machine Learning from Adaptive- and Multi-Grid
  {RANS} Simulations,} \emph{Proceedings of the 33rd Symposium on Naval
  Hydrodynamics, Osaka, Japan}, 2020{\natexlab{b}}.

\bibitem[{Serani et~al.(2016)Serani, Leotardi, Iemma, Campana, Fasano, and
  Diez}]{serani2016-ASC}
Serani, A., Leotardi, C., Iemma, U., Campana, E.~F., Fasano, G., and Diez, M.,
  \enquote{Parameter selection in synchronous and asynchronous deterministic
  particle swarm optimization for ship hydrodynamics problems,} \emph{Applied
  Soft Computing}, Vol.~49, 2016, pp. 313--334.

\bibitem[{Dykstra and Robertson(1982)}]{dykstra1982}
Dykstra, R.~L., and Robertson, T., \enquote{An Algorithm for Isotonic
  Regression for Two or More Independent Variables,} \emph{The Annals of
  Statistics}, Vol.~10, No.~3, 1982, pp. 708--716.

\end{thebibliography}

\end{document}